\documentclass[11pt, a4paper]{article}

\usepackage[utf8]{inputenc}
\usepackage[T1]{fontenc}
\usepackage{lmodern}
\usepackage{amssymb, amsthm, amsfonts, mathrsfs}
\usepackage{geometry}
\usepackage{mathtools}
\usepackage[all]{xy}
\usepackage{graphicx}
\usepackage{microtype}
\usepackage{enumitem}
\usepackage{stmaryrd}
\usepackage{amscd}
\usepackage[hidelinks]{hyperref}

\theoremstyle{plain}
\newtheorem{theorem}{Theorem}[section]
\newtheorem{lemma}[theorem]{Lemma}
\newtheorem{proposition}[theorem]{Proposition}
\newtheorem{corollary}[theorem]{Corollary}
\newtheorem*{proposition*}{Proposition}
\theoremstyle{definition}
\newtheorem{definition}[theorem]{Definition}
\newtheorem{example}[theorem]{Example}

\newtheorem{conjecture}[theorem]{Conjecture}
\theoremstyle{remark}
\newtheorem{remark}[theorem]{Remark}

\newtheorem*{theorem*}{Theorem}

\newcommand{\Z}{\mathbb{Z}}

\newcommand{\C}{\mathbb{C}}
\newcommand{\N}{\mathbb{N}}

\newcommand{\RP}{\mathbb{RP}}
\newcommand{\CP}{\mathbb{CP}}
\newcommand{\PP}{\mathbb{P}}

\newcommand{\Deck}{\operatorname{Deck}}
\newcommand{\Aut}{\operatorname{Aut}}
\newcommand{\Hom}{\operatorname{Hom}}
\newcommand{\colim}{\operatorname*{colim}}
\newcommand{\liminv}{\operatorname*{\varprojlim}}

\newcommand{\pihat}{\widehat{\pi}_1}
\newcommand{\cont}{\mathrm{cont}}
\newcommand{\W}{\mathcal{W}}
\newcommand{\PiST}{\Pi^{\mathrm{ST}}}
\newcommand{\HST}{H_{\mathrm{ST}}}

\title{Semi-topological Galois cohomology and Weierstrass realizability}
\author{Jyh-Haur Teh\thanks{Department of Mathematics, National Tsing Hua University, Taiwan.
E-mail: \texttt{jyhhaur@math.nthu.edu.tw}.}}
\date{}

\begin{document}
\maketitle

\begin{abstract}
Semi-topological Galois theory associates a canonical finite splitting covering to a monic Weierstrass polynomial. The inverse limit of the corresponding deck groups defines the absolute semi-topological Galois group, $\PiST(X,x)$. This paper develops a cohomology theory for $\PiST(X,x)$ with discrete torsion coefficients, establishing its fundamental properties and canonical comparison maps to singular cohomology. A Lyndon-Hochschild-Serre spectral sequence is used to yield an obstruction theory for semi-topological embedding problems. We prove several structural and vanishing results, including ST-fullness for free fundamental groups and triviality for finite fundamental groups. As applications, we provide a criterion for lifting finite projective monodromy to linear monodromy, formulate the $\pi_1$-detectable Weierstrass realizability conjecture for divisor classes and show that this conjecture is true for abelian varieties, smooth complex projective curves and ruled surfaces over positive-genus curves.
\end{abstract}

\section{Introduction}

The relationship between the topological classification of covering spaces and the algebraic data of polynomial equations is a central theme in geometry. Classical covering space theory organizes finite coverings of a space $X$ into a Galois category governed by the profinite completion $\pihat(X, x)$. Semi-topological Galois theory (\cite{LT15}, \cite{HansenBraidsCoverings}) introduces a powerful geometric refinement by distinguishing coverings that arise from Weierstrass polynomials.

For a Weierstrass polynomial $f \in C(X)[z]$ where $C(X)$ is the ring of complex-valued continuous functions on $X$, we can construct a canonical splitting covering $E_f \rightarrow X$, which is the minimal finite covering space where $f$ splits completely into distinct linear factors. The inverse limit of the deck groups $\Deck(E_f/X)$ over all splitting coverings forms a profinite group, the absolute semi-topological Galois group $\PiST(X,x)$. The canonical quotient map $\rho: \pihat(X, x) \rightarrow \PiST(X,x)$ precisely measures the discrepancy between the classical fundamental group and the semi-topological regime.

A motivating problem for this theory is the realization of geometric classes, such as divisors on projective manifolds, via polynomial data. To address this, we formulate the $\pi_1$-detectable Weierstrass realizability conjecture: every divisor class modulo $m$ detected by $\pi_1(X, x)$ lies in the image of the semi-topological comparison map.

To systematically study such realization problems, this paper develops a cohomology theory attached to $\PiST(X,x)$. For a discrete $\PiST(X,x)$-module $A$, we define the semi-topological Galois cohomology $H_{ST}^n(X,A)$ and construct canonical comparison maps $\Phi_X^n: H_{\mathrm{ST}}^n(X,A) \rightarrow H_{\mathrm{sing}}^n(X,A)$ to relate it to classical invariants.

We establish the fundamental properties of $H_{ST}^{\bullet}(X,A)$, including a continuous Lyndon-Hochschild-Serre spectral sequence that yields an obstruction theory for semi-topological embedding problems in degree 2. We then prove several structural theorems: ST-fullness for free fundamental groups, and the triviality of $\PiST(X, x)$ for spaces with finite fundamental groups. Crucially, for tori and abelian varieties, we prove that $\Phi_X^2$ is surjective for all finite coefficients, establishing our realizability conjecture for divisor classes on abelian varieties. We also show that the conjecture holds on smooth complex projective curves and ruled surfaces over positive-genus curves.

Beyond divisor realization, we show that semi-topological Galois cohomology controls the linearization of finite projective monodromy. We prove that finite projective semi-topological monodromy lifts to genuine linear monodromy precisely when the corresponding finite central-extension obstruction vanishes in semi-topological Galois cohomology.

\subsection*{Organization of the paper}
Section 2 reviews the semi-topological Galois theory of Weierstrass polynomials and splitting coverings, formulates the associated Galois category, and introduces the absolute semi-topological Galois group $\Pi^{ST}(X,x)$. It also defines the semi-topological Galois cohomology groups $H^{\bullet}_{ST}(X,A)$, establishes their basic properties, and constructs the comparison maps to singular cohomology together with the Lyndon--Hochschild--Serre spectral sequence attached to the kernel of $\widehat{\pi}_1(X,x)\to \Pi^{ST}(X,x)$. Section 3 studies semi-topological embedding problems, including low-degree interpretations of $H^{\bullet}_{ST}(X,A)$, braid-monodromy considerations, and the degree-two obstruction theory for central extensions. Section 4 is devoted to computations and structural results: ST-fullness for spaces with free fundamental group, vanishing results in the braid-free and torsion cases, the torus case, and an application to the linearization of finite projective monodromy. Finally, Section 5 formulates the $\pi_1$-detectable Weierstrass realizability conjecture and proves it for abelian varieties, smooth complex projective curves, and ruled surfaces over positive-genus curves.

\subsection*{Acknowledgments}
The author would like to thank the referee for a careful reading of the manuscript and for many helpful comments and suggestions, which have significantly improved the clarity and exposition of the paper. Any remaining errors or inaccuracies are the sole responsibility of the author. This work was supported by the National Science and Technology Council (Taiwan) under Grant No. NSTC 114-2115-M-007-002.

\section{The semi-topological Galois category and $\PiST(X,x)$}
\subsection{Weierstrass polynomials and splitting coverings}
We make a brief review of the semi-topological Galois theory developed in \cite{LT15}.
Throughout this paper, X denotes a topological space that is Hausdorff, path-connected, locally path-connected, and semilocally simply connected.
We write $C(X)$ for the ring of continuous functions $X\to\C$ and $C(X)[z]$ for the polynomial ring.

\begin{definition}[Weierstrass polynomial]
A polynomial
$$f(z)=a_nz^n+a_{n-1}z^{n-1}+\cdots+a_0\in C(X)[z]$$
is monic if its leading coefficient $a_n$ is the constant function $1_X:X\to \C, x\mapsto 1$.

A monic polynomial
\[
f(z)=z^n+a_{n-1}z^{n-1}+\cdots+a_0\in C(X)[z]
\]
is a \emph{Weierstrass polynomial of degree $n$ on $X$} if for every $x\in X$, the
specialization $f_x(z)\in\C[z]$ has $n$ distinct complex roots.
\end{definition}

\begin{definition}[Solution space and splitting in a covering]
Let $f\in C(X)[z]$ be a Weierstrass polynomial of degree $n$. The \emph{solution space} is the
subspace
\[
S_f:=\{(x,z)\in X\times \C\mid f_x(z)=0\}
\]
with projection $\pi:S_f\to X$, $\pi(x,z)=x$. If $p:Y\to X$ is a covering,
we say \emph{$f$ splits in $Y$} if $p^*f\in C(Y)[z]$ has $n$ distinct roots in $C(Y)$.
\end{definition}

\begin{definition}[Splitting covering]
A \emph{splitting covering} of $f$ is a path-connected covering $q:E_f\to X$
such that:
\begin{enumerate}[label=(\arabic*)]
\item $f$ splits in $E_f$,
\item $E_f$ is \emph{minimal} among such path-connected coverings: if $p:Y\to X$ is a covering
in which $f$ splits, then there is a (necessarily unique up to isomorphism over $X$)
covering map $Y\to E_f$ over $X$.
\end{enumerate}
\end{definition}

The minimality property is the topological analogue of the splitting field.
We will use the following fundamental facts as input.

\begin{theorem}[{\cite[Theorem 2.14]{LT15}}]\label{thm:split-exists}
For any Weierstrass polynomial $f$ on $X$, the following hold:
\begin{enumerate}[label=(\arabic*)]
\item A splitting covering $q:E_f\to X$ exists.
\item It is unique up to covering isomorphism over $X$.
\item The covering $E_f\to X$ is finite and Galois over $X$.
\end{enumerate}
\end{theorem}

\begin{definition}[Deck group]
For a covering $p:Y\to X$ we write
\[
\Deck(Y/X):=\{\varphi:Y\to Y\mid \varphi \text{ is a homeomorphism and }p\circ\varphi=p\}.
\]
\end{definition}

The cohomology theory we build depends on a profinite group extracted from splitting
coverings.

\subsection{The indexing system of splitting coverings}
Recall that a continuous map $p:Y\to X$ is called a \emph{covering map}
if every point $x\in X$ has an open neighborhood $U$ such that
\[
p^{-1}(U)=\coprod_{\lambda\in\Lambda} V_\lambda,
\]
where each $V_\lambda$ is open in $Y$ and the restriction
\[
p|_{V_\lambda}:V_\lambda\longrightarrow U
\]
is a homeomorphism. A covering map is called \emph{finite} if each
fiber $p^{-1}(x)$ is finite. Since throughout this paper $X$ is
path-connected, a finite covering of $X$ has a constant finite number
of sheets.

\begin{definition}
Fix a basepoint $x\in X$. Let $\W(X,x)$ be the category defined as follows:
\begin{itemize}
\item Objects are triples $(f,E_f,e_f)$
where $f$ is a Weierstrass polynomial on $X$, $E_f\to X$ is its splitting
covering, and $e_f\in E_f$ lies over $x$.

\item A morphism $(f,E_f,e_f)\longrightarrow (g,E_g,e_g)$ is a covering map $\phi:E_f\to E_g$ over $X$ such that
$\phi(e_f)=e_g$.
\end{itemize}
If such a morphism exists, we say that $(E_f,e_f)$
\emph{dominates} $(E_g,e_g)$. When the basepoint is understood, we shall
occasionally write $\W(X)$ for $\W(X,x)$.
\end{definition}

A morphism between connected covering spaces over $X$ is uniquely determined
by the image of a single point. Hence there is at most one morphism between
any two objects of $\W(X,x)$.

\begin{lemma}\label{lem:cofiltered}
The category $\W(X,x)$ is cofiltered. More precisely, for any two objects
$(f,E_f,e_f)$ and $(g,E_g,e_g)$ there exists an object $(h,E_h,e_h)$ and
morphisms such that
\[
(h,E_h,e_h)\longrightarrow(f,E_f,e_f),
\qquad
(h,E_h,e_h)\longrightarrow(g,E_g,e_g).
\]
\end{lemma}

\begin{proof}
Write
\[
f(z)=z^r+a_{r-1}z^{r-1}+\cdots+a_0,
\qquad
g(z)=z^s+b_{s-1}z^{s-1}+\cdots+b_0.
\]
For each $t\in X$, we set
\[
R_f(t):=2+\sum_{i=0}^{r-1}|a_i(t)|,
\qquad
R_g(t):=2+\sum_{j=0}^{s-1}|b_j(t)|,
\qquad
M(t):=R_f(t)+R_g(t)+1.
\]
By the Cauchy root bound(\cite[Chapter 8]{RahmanSchmeisser2002}), every root of $f_t$ lies in
$\overline{D(0,R_f(t))}$ and every root of $g_t$ lies in
$\overline{D(0,R_g(t))}$ where $\overline{D(0, r)}$ denotes the closed disc centered at the origin with radius $r$. Define
\[
h(z):=f(z)\,g(z-M)\in C(X)[z].
\]
For each $t\in X$, the roots contributed by $f_t$ lie in
$\overline{D(0,R_f(t))}$, while those contributed by $g_t(z-M(t))$
lie in $\overline{D(M(t),R_g(t))}$. Since
$M(t)>R_f(t)+R_g(t)$, these two closed disks are disjoint. Each factor
has only simple roots, so $h_t$ has $r+s$ distinct roots. Thus $h$ is a
Weierstrass polynomial on $X$.

Let $E_h\to X$ be the splitting covering of $h$. Since the two root clusters
above are disjoint and separated continuously, a splitting of $h$ on $E_h$
induces a splitting of $f$ and, after translating the second cluster by
$-M$, a splitting of $g$. By the minimality of splitting coverings, there
are covering maps
\[
\alpha:E_h\longrightarrow E_f,
\qquad
\beta:E_h\longrightarrow E_g
\]
over $X$.

Choose $e_h\in E_h$ over $x$. Because $E_f\to X$ and $E_g\to X$ are
connected Galois coverings, their deck groups act transitively on the fibres
over $x$. Hence, after postcomposing $\alpha$ and $\beta$ with suitable deck
transformations of $E_f$ and $E_g$, respectively, we may arrange
$\alpha(e_h)=e_f$ and $\beta(e_h)=e_g$. Thus these maps are morphisms in
$\W(X,x)$.

Finally, there are no nontrivial parallel pairs to equalize, because two
morphisms of connected pointed coverings over $X$ that agree at the chosen
point are equal. Therefore $\W(X,x)$ is cofiltered.
\end{proof}

\subsection{Definition of the absolute semi-topological Galois group}

For an object $(f,E_f,e_f)$ of $\W(X,x)$, set
\[
G_f:=\Deck(E_f/X).
\]
Suppose that
\[
\phi:(f,E_f,e_f)\longrightarrow(g,E_g,e_g)
\]
is a morphism. For the covering maps $q_f: E_f \rightarrow X$ and $q_g: E_g \rightarrow X$ , let
\[
H_f:=(q_f)_*\pi_1(E_f,e_f),
\qquad
H_g:=(q_g)_*\pi_1(E_g,e_g)
\]
be the corresponding subgroups of $\pi_1(X,x)$. Since $\phi$ is a pointed
map over $X$, we have $H_f\subseteq H_g$. By
Theorem~\ref{thm:split-exists}(3), both coverings are Galois, so $H_f$ and
$H_g$ are normal. Under the standard identifications
\[
G_f\cong \pi_1(X,x)/H_f,
\qquad
G_g\cong \pi_1(X,x)/H_g,
\]
the inclusion $H_f\subseteq H_g$ therefore induces a canonical surjective
homomorphism
\[
\mathrm{res}_{f,g}:G_f\twoheadrightarrow G_g.
\]
Equivalently, $\mathrm{res}_{f,g}(\sigma)$ is characterized by
\[
\phi\circ\sigma
=
\mathrm{res}_{f,g}(\sigma)\circ\phi
\qquad(\sigma\in G_f).
\]
Thus the deck groups of pointed splitting coverings form an inverse system of
finite groups over the cofiltered category $\W(X,x)$.

\begin{definition}[Absolute semi-topological Galois group]\label{def:PiST}
The \emph{absolute semi-topological Galois group} of $X$ based at $x$ is
\[
\PiST(X,x)
:=
\liminv_{(f,E_f,e_f)\in\W(X,x)}\Deck(E_f/X),
\]
endowed with its inverse-limit topology. In particular, $\PiST(X,x)$ is a
profinite group.
\end{definition}

Our definition of the absolute semi-topological Galois group is parallel to the classical absolute Galois group of a field.
In this analogy, a splitting covering associated with a Weierstrass polynomial plays the role of a splitting field of a polynomial.
Thus, finite coverings dominated by splitting coverings are the topological counterparts of finite Galois extensions,
and their inverse limit naturally gives the absolute semi-topological Galois group.

\begin{remark}[Dependence on the basepoint]
If $x'$ is another basepoint and $\gamma$ is a path from $x$ to $x'$, lifting
$\gamma$ in each splitting covering identifies the corresponding pointed
inverse systems and hence gives an isomorphism
$\PiST(X,x)\cong\PiST(X,x')$. Changing the path $\gamma$ changes this
isomorphism by an inner automorphism. Thus, as in ordinary fundamental-group
theory, the group is independent of the basepoint up to inner automorphism.
\end{remark}

\subsection{Galois category formulation}

Let $\operatorname{Cov}_{\mathrm{fin}}(X)$ denote the category of finite
covering spaces of $X$, not necessarily connected, with morphisms given by
maps over $X$. Define $\operatorname{Cov}_{\mathrm{ST}}(X)$ to be the full
subcategory consisting of those finite coverings $Y\to X$ for which every
connected component of $Y$ is dominated by a splitting covering. Equivalently,
$Y\to X$ lies in $\operatorname{Cov}_{\mathrm{ST}}(X)$ if there exist
Weierstrass polynomials $f_1,\ldots,f_m$ and a surjective map of coverings
\[
\coprod_{i=1}^m E_{f_i}\longrightarrow Y
\]
over $X$.

It is useful to distinguish $\W(X,x)$ from $\operatorname{Cov}_{\mathrm{ST}}(X)$.
The category $\W(X,x)$ consists of connected pointed Galois coverings and serves
as the cofiltered indexing category in the definition of $\PiST(X,x)$,
whereas $\operatorname{Cov}_{\mathrm{ST}}(X)$ denotes the category of all finite
coverings, not necessarily connected, dominated by objects of $\W(X,x)$.

The basepoint $x$ defines the fibre functor
\[
\Phi_x:\operatorname{Cov}_{\mathrm{ST}}(X)\longrightarrow
\operatorname{FinSet},
\qquad
(Y\xrightarrow{p}X)\longmapsto p^{-1}(x).
\]

Recall that a \emph{Galois category} is a category $\mathcal{C}$ together
with a fibre functor
\[
F:\mathcal{C}\longrightarrow \operatorname{FinSet},
\]
satisfying the standard Grothendieck axioms: the category $\mathcal{C}$
admits finite limits and finite coproducts, every morphism admits a
factorization as a strict epimorphism followed by a monomorphism, and
the functor $F$ preserves the relevant finite limits and coproducts and
is conservative. See \cite{SGA1} or \cite{BorceuxJanelidzeGalois} for the complete list of
axioms.

\begin{theorem}\label{thm:galois-cat}
The pair $(\operatorname{Cov}_{\mathrm{ST}}(X),\Phi_x)$ is a Galois
category. Moreover, there is a natural isomorphism of profinite groups
\[
\Aut(\Phi_x)\cong\PiST(X,x),
\]
and $\operatorname{Cov}_{\mathrm{ST}}(X)$ is equivalent, compatibly with
$\Phi_x$, to the category of finite sets endowed with a continuous
$\PiST(X,x)$-action.
\end{theorem}

\begin{proof}
For each object $(f,E_f,e_f)$ of $\W(X,x)$, the
projection from the inverse limit gives a continuous surjection
$
\PiST(X,x)\twoheadrightarrow G_f.
$
If a connected finite covering $Y\to X$ is dominated by $E_f$, then the
monodromy action of $\pi_1(X,x)$ on the fibre $Y_x$ factors through the finite
quotient $G_f$. Hence $Y_x$ acquires a continuous $\PiST(X,x)$-action. For a general
object of $\operatorname{Cov}_{\mathrm{ST}}(X)$, Lemma~\ref{lem:cofiltered}
(and induction on the number of connected components) provides a single
splitting covering dominating all its components. Consequently the fibre
functor $\Phi_x$ lifts to a functor
\[
\operatorname{Cov}_{\mathrm{ST}}(X)
\longrightarrow \PiST(X,x)\text{-}\operatorname{FinSet}.
\]

Conversely, let $S$ be a finite continuous $\PiST(X,x)$-set. Its action factors through
a finite quotient of $\PiST(X,x)$. Since $\PiST(X,x)$ is an inverse limit of finite groups with
surjective transition maps, the kernels of the projections $\PiST(X,x)\to G_f$ form a
neighbourhood basis of the identity. Thus the action on $S$ factors through
$G_f$ for some $(f,E_f,e_f)$. By the classical classification of finite
coverings by finite $\pi_1(X,x)$-sets, $S$ therefore corresponds to a finite
covering of $X$ dominated by $E_f$, and hence to an object of
$\operatorname{Cov}_{\mathrm{ST}}(X)$.

These two constructions are mutually quasi-inverse and identify $\Phi_x$ with
the forgetful functor $\PiST(X,x)\text{-}\operatorname{FinSet}\to\operatorname{FinSet}$.
The category of finite continuous $\PiST(X,x)$-sets is a Galois category, and the
automorphism group of its forgetful fibre functor is naturally $\PiST(X,x)$. The
assertions follow.
\end{proof}

\subsection{Continuous cochains and continuous cohomology}

Let $G$ be a profinite group. A \emph{discrete $G$-module} is an abelian
group $A$, endowed with the discrete topology and a $G$-action, such that the
action map $G\times A\to A$ is continuous. Equivalently, the stabilizer of
every $a\in A$ is open in $G$; equivalently,
\[
A=\bigcup_{U\triangleleft_o G}A^U,
\qquad
A^U:=\{a\in A:ua=a\text{ for every }u\in U\},
\]
where $U\triangleleft_o G$ means that $U$ is an open normal subgroup of $G$.

For $n\geq0$, let
\[
C^n_{\mathrm{cont}}(G,A):=\operatorname{Map}_{\mathrm{cont}}(G^n,A),
\]
with $C^0_{\mathrm{cont}}(G,A)=A$, and equip these groups with the usual
inhomogeneous differential
\[
\begin{aligned}
(dc)(g_1,\ldots,g_{n+1})={}&g_1c(g_2,\ldots,g_{n+1})+\sum_{i=1}^{n}(-1)^ic(g_1,\ldots,g_ig_{i+1},\ldots,g_{n+1})
+(-1)^{n+1}c(g_1,\ldots,g_n)
\end{aligned}
\]

\begin{definition}[Continuous cohomology]\label{def:cont-coh}
The \emph{continuous cohomology} of $G$ with coefficients in the discrete
$G$-module $A$ is
\[
H^n_{\mathrm{cont}}(G,A):=H^n\bigl(C^\bullet_{\mathrm{cont}}(G,A)\bigr)
\]
\end{definition}

The following finite-quotient description is standard; see, for example,
\cite[Chapter~6]{RibesZalesskii} or \cite{SerreGaloisCohom}.

\begin{proposition}\label{prop:finite-quot-formula}
Let $G$ be a profinite group and $A$ a discrete $G$-module. Then, for every
$n\geq0$,
\[
H^n_{\mathrm{cont}}(G,A)\cong \varinjlim_{U\triangleleft_o G} H^n(G/U,A^U).
\]
If $V\subseteq U$ are open normal subgroups, the transition map
\[
H^n(G/U,A^U)\longrightarrow H^n(G/V,A^V)
\]
is induced by inflation along $G/V\twoheadrightarrow G/U$ together with the
inclusion $A^U\hookrightarrow A^V$. In particular, if the $G$-action on $A$
is trivial, then
\[
H^n_{\mathrm{cont}}(G,A)
\cong
\varinjlim_{U\triangleleft_o G}H^n(G/U,A).
\]
\end{proposition}

\begin{proof}
Every continuous map from the compact space $G^n$ to the discrete space $A$
has finite image and is locally constant. Hence each continuous cochain
factors through $(G/U)^n$ for some open normal subgroup $U$ that fixes all
values of the cochain. It follows that
\[
C^\bullet_{\mathrm{cont}}(G,A)
\cong
\varinjlim_{U\triangleleft_o G}C^\bullet(G/U,A^U)
\]
as cochain complexes. Filtered colimits are exact in the category of abelian
groups, so taking cohomology commutes with this colimit and gives the stated
isomorphism.
\end{proof}

\subsection{Definition of $\HST^\bullet(X,A)$}

For each object $(f,E_f,e_f)$ of $\W(X,x)$, let
\[
\pi_f:\PiST(X,x)\longrightarrow \Deck(E_f/X)
\]
be the canonical projection. Since the transition homomorphisms between the
deck groups are surjective, each $\pi_f$ is surjective. Hence
\[
\PiST(X,x)/\ker(\pi_f)\cong \Deck(E_f/X).
\]

\begin{definition}[Semi-topological Galois cohomology]\label{def:HST}
Let $A$ be a discrete $\PiST(X,x)$-module. The \emph{semi-topological
Galois cohomology} of $X$ with coefficients in $A$ is defined by
\[
\HST^n(X,A):=H^n_{\mathrm{cont}}(\PiST(X,x),A)
\]
\end{definition}

The preceding definition admits a concrete description in terms of the
finite deck groups associated with splitting coverings. We record this
description explicitly, since the direction of the transition maps will be
important later.

\begin{proposition}\label{prop:HST-finite-stage}
Let $A$ be a discrete torsion $\PiST(X,x)$-module. Then there is a natural
isomorphism
\[
\HST^n(X,A)\cong\colim_{(f,E_f,e_f)\in\W(X,x)^{\mathrm{op}}}
H^n\bigl(\Deck(E_f/X),A^{\ker(\pi_f)}\bigr),
\]
where
\[
A^{\ker(\pi_f)}
:=
\left\{a\in A \ |\ ka=a \text{ for every } k\in\ker(\pi_f)
\right\}.
\] and $\W(X,x)^{\mathrm{op}}$ is the opposite category of $\W(X, x)$.
\end{proposition}

This is the standard continuity property of the cohomology of profinite groups; see \cite[Proposition 1.2.5]{NeukirchSchmidtWingberg}.

\begin{remark}[Direction of the transition maps]
The appearance of $\W(X,x)^{\mathrm{op}}$ is essential. A morphism
\[
(f,E_f,e_f)\longrightarrow(g,E_g,e_g)
\]
in $\W(X,x)$ induces a surjective homomorphism
\[
\Deck(E_f/X)\longrightarrow\Deck(E_g/X).
\]
Since group cohomology is contravariant with respect to group
homomorphisms, the corresponding homomorphism on cohomology goes in the
opposite direction:
\[
H^n\bigl(\Deck(E_g/X),A\bigr)
\longrightarrow
H^n\bigl(\Deck(E_f/X),A\bigr).
\]
Thus the transition homomorphisms in the finite-stage description of
$\HST^\bullet(X,A)$ are inflation homomorphisms, rather than restriction
homomorphisms.
\end{remark}

\begin{remark}[Sheaf-theoretic viewpoint]
Endow $\mathrm{Cov}_{\mathrm{ST}}(X)$ with its canonical Galois topology.
By Theorem~\ref{thm:galois-cat}, the associated Galois category is governed
by the profinite group $\PiST(X,x)$. In particular, finite objects correspond
to finite sets equipped with a continuous $\PiST(X,x)$-action. Passing to
the associated classifying topos identifies discrete $\PiST(X,x)$-modules
with the corresponding discrete abelian sheaves. Under this identification,
\[
\HST^\bullet(X,A)=H^\bullet_{\mathrm{cont}}(\PiST(X,x),A)
\]
is the cohomology of the associated classifying topos with coefficients in
the sheaf determined by $A$.
\end{remark}

\subsection{Properties}

\begin{proposition}\label{prop:functoriality}
Let $X$ and $Y$ be Hausdorff, path-connected, locally path-connected and semilocally simply connected topological spaces.
Suppose $u:(Y, y)\rightarrow (X, x)$ is a continuous map preserving basepoints.
Then pullback of coverings defines a functor
$
u^*:\mathrm{Cov}_{\mathrm{ST}}(X)\to \mathrm{Cov}_{\mathrm{ST}}(Y),
$
hence a continuous homomorphism of profinite groups
\[
u_\#: \PiST(Y,y)\longrightarrow \PiST(X,x)
\]
Consequently, for any discrete $\PiST(X, x)$-module $A$ there are natural pullback maps
\[
u^*:\HST^n(X,A)\to \HST^n(Y,u^*A)
\]
Here $u^*A$ denotes $A$ regarded as a discrete $\PiST(Y,y)$-module via $u_{\#}$; explicitly,
\[
g\cdot a:=u_\#(g)\cdot a
\]
for $g\in\PiST(Y,y)$ and $a\in A$.
\end{proposition}

\begin{proof}[Proof sketch]
If $E_f\to X$ is a splitting covering, its pullback $u^*E_f\to Y$ is a finite cover of $Y$ which
belongs to $\mathrm{Cov}_{\mathrm{ST}}(Y)$. The induced homomorphism on automorphism groups of fibre functors gives $u_\#$.
Functoriality on cohomology follows from functoriality of continuous cohomology.
\end{proof}

\begin{proposition}\label{prop:les}
A short exact sequence of discrete $\PiST(X, x)$-modules
$0\to A'\to A\to A''\to 0$
induces a natural long exact sequence
\[
\cdots\to \HST^n(X,A')\to \HST^n(X,A)\to \HST^n(X,A'')\to \HST^{n+1}(X,A')\to\cdots.
\]
\end{proposition}
This follows directly from the standard long exact sequence for continuous cohomology of profinite groups \cite[Chapter 6, §6.2]{RibesZalesskii}.

\begin{proposition}\label{prop:cup}
For discrete $\PiST(X, x)$-modules $A,B$ there are natural cup products
\[
\smile:\HST^p(X,A)\times \HST^q(X,B)\to \HST^{p+q}(X,A\otimes B),
\]
which are graded-commutative when $A=B$ and the action is trivial.
\end{proposition}

For the general definition and properties of cup products in profinite group cohomology, see \cite[Chapter 6, §6.7]{RibesZalesskii}.

\subsection{Comparison with classical topology}
Connected Galois finite covers of $X$ are classified by finite quotients of $\pi_1(X,x)$, hence by
the profinite completion $\pihat(X,x)$.
Since the splitting-cover quotients form a cofiltered inverse system: the image of $\pihat(X, x)$ in $\PiST(X,x)$
under the canonical homomorphism is compact, hence closed, and its projection to every finite stage is surjective;
cofilteredness makes the image dense in the inverse limit, hence equal to the whole inverse limit.

\begin{proposition}\label{prop:rho}
There is a natural continuous surjective homomorphism of profinite groups
\[
\rho:\pihat(X,x)\twoheadrightarrow \PiST(X,x),
\]
characterized by the property that for each splitting covering $E_f\to X$ the composite
\[
\pihat(X,x)\xrightarrow{\rho}\PiST(X,x)\to \Deck(E_f/X)
\]
is the usual surjection corresponding to the finite quotient defined by $E_f\to X$.
\end{proposition}

\begin{definition}[Semi-topological fundamental kernel]
Define the closed normal subgroup
\[
K_{\mathrm{ST}}(X,x):=\ker(\rho)\subset \pihat(X,x).
\]
Equivalently, $K_{\mathrm{ST}}(X,x)$ is the intersection of open normal subgroups of
$\pihat(X,x)$ corresponding to splitting coverings.
\end{definition}

\subsection{A comparison map to singular cohomology}

Let $A$ be a finite abelian group with trivial action.
The morphism $\rho$ induces a pullback map
$
\rho^*:H^n_{\mathrm{cont}}(\PiST(X, x),A)\to H^n_{\mathrm{cont}}(\pihat(X, x),A).
$
The canonical map $\iota:\pi_1(X,x)\to \pihat(X,x)$ induces
$
\iota^*:H^n_{\mathrm{cont}}(\pihat(X, x),A)\to H^n(\pi_1(X,x),A).
$
Finally, the classifying map $c:X\to B\pi_1(X,x)$ gives
$c^*:H^n(\pi_1(X,x),A)\to H^n_{\mathrm{sing}}(X,A)$. Note that the singular cohomology
of the classifying space $B\pi_1(X,x)$ is the group cohomology of $\pi_1(X,x)$.

\begin{definition}[Comparison map and realizability]\label{def:comparison}
For finite abelian group $A$ with trivial action define
\[
\Phi^n_X:\HST^n(X,A)\xrightarrow{\rho^*}H^n_{\mathrm{cont}}(\pihat(X, x),A)
\xrightarrow{\iota^*}H^n(\pi_1(X, x),A)\xrightarrow{c^*}H^n_{\mathrm{sing}}(X,A).
\]
A class $\alpha\in H^n_{\mathrm{sing}}(X,A)$ is \emph{Weierstrass realizable} if
$\alpha\in \mathrm{Im}(\Phi^n_X)$.
\end{definition}

\begin{remark}[Aspherical case]
Recall that a path-connected topological space $X$ is called
\emph{aspherical} if $\pi_n(X)=0$ for all $n\geq 2$.
Equivalently, if $X$ has the homotopy type of a CW-complex, then $X$
is aspherical if and only if its universal covering space is
contractible. In this case, $X$ has the homotopy type of an
Eilenberg--Mac Lane space $K(\pi_1(X, x),1)$.
If $X$ is a $K(\pi,1)$, then $c^*$ is an isomorphism and $\Phi^n_X$ lands in
$H^n(\pi_1(X, x),A)\cong H^n_{\mathrm{sing}}(X,A)$ without ambiguity.
\end{remark}

\begin{theorem}[Naturality of the comparison map]\label{thm:comparison-naturality}
Let
$
u:(Y,y)\longrightarrow (X,x)
$
be a continuous based map, and let $A$ be a finite abelian group with trivial action. Then the diagram
\[
\xymatrix{
\HST^n(X,A) \ar[r]^{\Phi_X^n} \ar[d]_{u^*} &
H^n_{\mathrm{sing}}(X,A) \ar[d]^{u^*} \\
\HST^n(Y,A) \ar[r]^{\Phi_Y^n} &
H^n_{\mathrm{sing}}(Y,A)
}
\]
commutes.
\end{theorem}

\begin{proof}
The result follows from the following two commutative diagrams:
\[
\xymatrix{
\pi_1(Y,y) \ar[r]^{\iota_Y} \ar[d]_{u_*} &
\widehat{\pi}_1(Y,y) \ar[r]^{\rho_Y} \ar[d]_{\widehat{u_*}} &
\PiST(Y,y) \ar[d]^{u_\#} \\
\pi_1(X,x) \ar[r]^{\iota_X} &
\widehat{\pi}_1(X,x) \ar[r]^{\rho_X} &
\PiST(X,x)
}
\]
and

\[
\xymatrix{
H^n_{\mathrm{cont}}(\PiST(X,x),A) \ar[r]^{\rho_X^*} \ar[d]_{u_\#^*} &
H^n_{\mathrm{cont}}(\widehat{\pi}_1(X,x),A) \ar[r]^{\iota_X^*} \ar[d]_{\widehat{u_*}^{\,*}} &
H^n(\pi_1(X,x),A) \ar[d]^{u_*^*} \\
H^n_{\mathrm{cont}}(\PiST(Y,y),A) \ar[r]^{\rho_Y^*} &
H^n_{\mathrm{cont}}(\widehat{\pi}_1(Y,y),A) \ar[r]^{\iota_Y^*} &
H^n(\pi_1(Y,y),A).
}
\]

\end{proof}

The exact sequence
$1\to K_{\mathrm{ST}}(X, x)\to \pihat(X, x)\to \PiST(X, x)\to 1$
yields the standard continuous Lyndon--Hochschild--Serre (LHS) spectral sequence.

\begin{theorem}\label{thm:lhs}
For any finite discrete $\PiST(X, x)$-module $A$ (with trivial action, for simplicity) there is a
first quadrant spectral sequence
\[
E_2^{p,q}=H^p_{\mathrm{cont}}(\PiST(X, x), H^q_{\mathrm{cont}}(K_{\mathrm{ST}}(X, x),A))
\Longrightarrow H^{p+q}_{\mathrm{cont}}(\pihat(X, x),A).
\]
In particular, the edge map $E_2^{n,0}\to H^n_{\mathrm{cont}}(\pihat(X, x),A)$ is $\rho^*$.
\end{theorem}

For this result, we refer the reader to \cite[Chapter 6, Theorem 6.8.2]{RibesZalesskii}.

\begin{remark}
The spectral sequence measures, cohomologically, the failure of $\rho$ to be an
isomorphism. When $K_{\mathrm{ST}}(X, x)$ is nontrivial, classes in
$H^q_{\mathrm{cont}}(K_{\mathrm{ST}}(X, x),A)$ can obstruct surjectivity of $\rho^*$.
\end{remark}

\section{Semi-topological embedding problems}

\subsection{A realizability criterion}

\begin{proposition}\label{prop:H0H1}
Let $A$ be a finite abelian group with trivial action.
\begin{enumerate}[label=(\arabic*)]
\item $\HST^0(X,A)\cong A$.
\item $\HST^1(X,A)\cong \Hom_{\mathrm{cont}}(\PiST(X, x),A)$.
\item If $A$ is viewed as a constant sheaf on $\mathrm{Cov}_{\mathrm{ST}}(X)$, then
$\HST^1(X,A)$ classifies $A$-torsors in $\mathrm{Cov}_{\mathrm{ST}}(X)$. The connected $A$-torsors
correspond precisely to surjective continuous homomorphisms $\PiST(X,x)\twoheadrightarrow A$,
and hence to connected Galois coverings with deck group $A$.
\end{enumerate}
\end{proposition}

\begin{proof}
Statements (1) and (2) are the standard descriptions of continuous
group cohomology in degrees zero and one. Since the action of
$\Pi^{\mathrm{ST}}(X, x)$ on $A$ is trivial, we have
\[
H^0_{\mathrm{ST}}(X,A)
=
H^0_{\mathrm{cont}}(\Pi^{\mathrm{ST}}(X, x),A)
=
A
\]
and
\[
H^1_{\mathrm{ST}}(X,A)
=
H^1_{\mathrm{cont}}(\Pi^{\mathrm{ST}}(X, x),A)
\cong
\operatorname{Hom}_{\mathrm{cont}}(\Pi^{\mathrm{ST}}(X, x),A).
\]
See \cite{SerreGaloisCohom} for these standard facts concerning continuous
cohomology of profinite groups.

For (3), Theorem~2.9 identifies
$\operatorname{Cov}_{\mathrm{ST}}(X)$ with the category of finite sets
equipped with a continuous action of $\Pi^{\mathrm{ST}}(X, x)$.
Under this equivalence, $A$-torsors in
$\operatorname{Cov}_{\mathrm{ST}}(X)$ correspond to continuous
$A$-torsors for $\Pi^{\mathrm{ST}}(X, x)$, whose isomorphism classes are
classified by
\[
H^1_{\mathrm{cont}}(\Pi^{\mathrm{ST}}(X, x),A)
=
H^1_{\mathrm{ST}}(X,A).
\]
The connected $A$-torsors are precisely the connected Galois coverings with deck group $A$.
\end{proof}

\begin{proposition}\label{prop:deg1-criterion}
Let $A$ be a finite abelian group with trivial action and let
$\chi\in H^1_{\mathrm{sing}}(X,A)$.
Assume $X$ is aspherical so that $H^1_{\mathrm{sing}}(X,A)\cong \Hom(\pi_1(X, x),A)$.
Then $\chi$ is Weierstrass realizable if and only if there exists a Weierstrass polynomial
$f$ such that $\ker(\chi)$ contains the subgroup of $\pi_1(X, x)$ corresponding to the
splitting covering $E_f\to X$ (equivalently, $\chi$ factors through $\Deck(E_f/X)$).
\end{proposition}

\begin{proof}
By Definition~\ref{def:HST}, a class in $\HST^1(X,A)$ is a continuous homomorphism
$\PiST(X, x)\to A$ and therefore factors through some finite quotient $\Deck(E_f/X)$.
Tracing through $\Phi^1_X$ yields exactly the factorization statement.
\end{proof}

\subsection{Degree $2$ and embedding problems}

For a finite group $G$ and a finite abelian group $A$ with trivial $G$-action,
the group $H^2(G,A)$ classifies equivalence classes of central extensions
$1\to A\to H\to G\to 1$.

\begin{definition}[Semi-topological embedding problem]
Fix a splitting covering $F\to X$ with group $G:=\Deck(F/X)$ and a central extension
$1\to A\to H\to G\to 1$.
The \emph{semi-topological embedding problem} asks whether there exists a splitting
covering $E\to X$ dominating $F$ such that $\Deck(E/X)\cong H$ and the induced
restriction $\Deck(E/X)\to \Deck(F/X)$ identifies with $H\twoheadrightarrow G$.

We say that the semi-topological embedding problem is \emph{solvable}
if such a splitting covering $E\to X$ exists. In this case,
$E\to X$, together with the above identification of its deck group
with $H$, is called a \emph{solution} of the embedding problem.
\end{definition}

The following is a semi-topological analogue of the standard cohomological obstruction for central embedding problems; 
see, for example, \cite[Proposition VI.16.1 ]{BerhuyGaloisCohomology}.

\begin{theorem}\label{thm:obstruction}
Let $A$ be a finite abelian group with trivial $G$-action, and let
\[
1\longrightarrow A\longrightarrow H
\stackrel{\pi}{\longrightarrow}G\longrightarrow1
\]
be a central extension with extension class
$[H]\in H^2(G,A)$.
Let $q:\Pi^{\mathrm{ST}}(X,x)\twoheadrightarrow G$
be the quotient homomorphism corresponding to the splitting covering
$F\to X$.

If the semi-topological embedding problem associated to
$F\to X$ and the above central extension is solvable, then
\[
q^*[H]=0
\qquad\text{in}\qquad
H^2_{\mathrm{ST}}(X,A),
\]
where
$q^*:H^2(G,A)\longrightarrow H^2_{\mathrm{cont}}(\Pi^{\mathrm{ST}}(X,x),A)= H^2_{\mathrm{ST}}(X,A)$
is the inflation map induced by $q$.
\end{theorem}

\begin{proof}
Assume that the semi-topological embedding problem is solvable.
By Definition~3.3, there exists a splitting covering
$E\longrightarrow X $
dominating $F\to X$, together with an identification
$
\operatorname{Deck}(E/X)\cong H
$
such that the induced restriction homomorphism
\[
\operatorname{Deck}(E/X)\longrightarrow
\operatorname{Deck}(F/X)
\]
identifies with the given surjection
$
\pi:H\twoheadrightarrow G.
$

Since $E\to X$ is a splitting covering, the definition of
$\Pi^{\mathrm{ST}}(X,x)$ gives a canonical continuous surjection
\[
\psi:
\Pi^{\mathrm{ST}}(X,x)
\twoheadrightarrow
\operatorname{Deck}(E/X)
\cong H
\]
Likewise, the splitting covering $F\to X$ determines the quotient
\[
q:
\Pi^{\mathrm{ST}}(X,x)
\twoheadrightarrow
\operatorname{Deck}(F/X)
\cong G.
\]
Because $E$ dominates $F$, these quotient maps are compatible, and
hence
$
\pi\circ\psi=q
$

Now pull back the central extension
\[
1\longrightarrow A\longrightarrow H
\stackrel{\pi}{\longrightarrow}G\longrightarrow1
\]
along $q$. This gives the central extension
\[
1\longrightarrow A\longrightarrow
\Pi^{\mathrm{ST}}(X,x)\times_G H
\longrightarrow
\Pi^{\mathrm{ST}}(X,x)
\longrightarrow1,
\]
whose cohomology class is
\[
q^*[H]\in
H^2_{\mathrm{cont}}(\Pi^{\mathrm{ST}}(X,x),A)
=
H^2_{\mathrm{ST}}(X,A).
\]

The homomorphism $\psi$ defines a continuous section
\[
s:\Pi^{\mathrm{ST}}(X,x)
\longrightarrow
\Pi^{\mathrm{ST}}(X,x)\times_G H,
\qquad
s(\gamma)=(\gamma,\psi(\gamma)),
\]
because
\[
q(\gamma)=\pi(\psi(\gamma)).
\]
Hence the pulled-back extension splits. By the classification of
central extensions by second continuous group cohomology, its
extension class is zero. Therefore
$
q^*[H]=0
$
in $H^2_{\mathrm{ST}}(X,A)$.
\end{proof}

\section{Computations and structural results}

\subsection{The free fundamental group case}

The semi-topological theory is particularly transparent when every finite cover is
dominated by a splitting covering.

\begin{definition}[ST-fullness]
We say that $X$ is \emph{ST-full} if every finite covering
$Y\to X$ belongs to $\operatorname{Cov}_{\mathrm{ST}}(X)$; equivalently,
every finite covering of $X$ is dominated by a finite disjoint union
of splitting coverings arising from Weierstrass polynomials on $X$.

Since $\operatorname{Cov}_{\mathrm{ST}}(X)$ is a full
subcategory of the category $\operatorname{Cov}_{\mathrm{fin}}(X)$ of
finite coverings of $X$, this condition is equivalent to
\[
\operatorname{Cov}_{\mathrm{ST}}(X)
=
\operatorname{Cov}_{\mathrm{fin}}(X).
\]
\end{definition}

\begin{theorem}\label{thm:free-stfull}
If $\pi_1(X, x)$ is a free group, then $X$ is ST-full. In particular,
\[
\PiST(X, x)\cong \pihat(X, x).
\]
\end{theorem}

\begin{proof}[Proof sketch]
By \cite[Theorem 6.3]{HansenBraidsCoverings}, if
$\pi_1(X, x)$ is free, every finite covering space of $X$ is equivalent to a solution space of
some Weierstrass polynomial on $X$.
Each connected component of the solution-space covering is dominated by the corresponding connected splitting covering; consequently the whole finite covering is dominated by a finite disjoint union of splitting coverings and therefore belongs to
$\mathrm{Cov}_{\mathrm{ST}}(X)$. Thus $\PiST(X, x)$
captures all finite covers and equals $\pihat(X, x)$.
\end{proof}

\begin{corollary}\label{cor:free-vanish}
Assume $\pi_1(X, x)$ is free (possibly of countable rank). For any finite abelian group $A$
with trivial action,
\[
\HST^n(X,A)=0\qquad (n\ge 2).
\]
\end{corollary}

\begin{proof}
By Theorem~\ref{thm:free-stfull}, $\PiST(X, x)\cong \pihat(X, x)$ is a free profinite group.
Free profinite groups have cohomological dimension $1$, hence $\HST^n(X,A)=0$
for $n\ge 2$.
\end{proof}

\subsection{Braid monodromy}
Let
\[
f(z)=z^n+a_{n-1}z^{n-1}+\cdots+a_0
\in C(X)[z]
\]
be a Weierstrass polynomial of degree $n$. For each $x\in X$, we denote
by
\[
f_x(z)
=
z^n+a_{n-1}(x)z^{n-1}+\cdots+a_0(x)
\in \mathbb{C}[z]
\]
the specialization of $f$ at $x$. Since $f$ is a Weierstrass polynomial,
$f_x$ has $n$ distinct complex roots for every $x\in X$.

Let $P_n^{\mathrm{simp}}$ denote the space of monic complex polynomials
of degree $n$ having distinct roots. The assignment
\[
a_f:X\longrightarrow P_n^{\mathrm{simp}},
\qquad
x\longmapsto f_x,
\]
defines a continuous map.

Fix a basepoint $x_0\in X$, and take $f_{x_0}$ as the corresponding
basepoint of $P_n^{\mathrm{simp}}$. The induced homomorphism
\[
\theta_f:=(a_f)_*:
\pi_1(X,x_0)
\longrightarrow
\pi_1(P_n^{\mathrm{simp}},f_{x_0})
\cong B_n
\]
is called the \emph{braid monodromy} of $f$, where $B_n$ denotes the
braid group on $n$ strands. Changing the basepoint changes this
homomorphism only up to the usual conjugacy arising from a change of
basepoint.

We will use the standard fact:

\begin{proposition}\label{prop:Bn-torsionfree}
For every $n\ge 1$, the braid group $B_n$ is torsion-free. In particular, any
homomorphism from a torsion group into $B_n$ is trivial.
\end{proposition}

\begin{definition}[Braid-free]\label{def:braidfree}
We call a group $\Gamma$ \emph{braid-free} if for every $n\ge 1$ every homomorphism
$\Gamma\to B_n$ is trivial. We call a space $X$ \emph{braid-free} if $\pi_1(X, x)$ is braid-free.
\end{definition}

\begin{remark}
If $\pi_1(X, x)$ is a torsion group, then $X$ is braid-free
by Proposition~\ref{prop:Bn-torsionfree}. This includes the case where $\pi_1(X, x)$ is finite.
\end{remark}

\begin{theorem}\label{thm:braidfree-trivial}
Let $X$ be a Hausdorff, path-connected, locally path-connected and semilocally simply connected topological space.
Assume that $X$ is braid-free. Then:
\begin{enumerate}[label=(\arabic*)]
\item Every splitting covering $E_f\to X$ is (isomorphic to) the identity covering.
\item Hence $\PiST(X,x)=1$ for every basepoint $x$.
\item Consequently, for every discrete coefficient module $A$ and every $n>0$,
\[
\HST^n(X,A)=H^n_{\cont}(\PiST(X, x),A)=0.
\]
\end{enumerate}
\end{theorem}

\begin{proof}
Fix a Weierstrass polynomial $f\in C(X)[z]$ of degree $n$.
Since $X$ is braid-free, the braid monodromy
$\theta_f:\pi_1(X,x)\to B_n$ must be the trivial homomorphism.

The standard surjective homomorphism
\[
\sigma:B_n\longrightarrow S_n
\]
records the permutation of the $n$ roots determined by a braid. Hence
the monodromy representation of the solution-space covering
\[
p:S_f\longrightarrow X
\]
is the composite
\[
\pi_1(X,x)
\xrightarrow{\;\theta_f\;}
B_n
\xrightarrow{\;\sigma\;}
S_n.
\]
Since each fiber $p^{-1}(x)$ consists of the $n$ distinct roots of
$f_x$, the map $p:S_f\to X$ is an $n$-sheeted covering.

If $\theta_f$ is trivial, then the associated permutation monodromy
$\sigma\circ\theta_f$ is also trivial. Therefore no permutation of the
$n$ sheets occurs under continuation along loops in $X$, and the
covering is trivial:
\[
S_f\cong \coprod_{i=1}^n X
\]
as a covering over $X$.
Equivalently, the $n$ roots of $f$ can be chosen globally as continuous
functions $r_1,\ldots,r_n\in C(X)$, so that
\[
f(z)=\prod_{i=1}^n (z-r_i)
\qquad\text{in } C(X)[z].
\]

Let $q:E_f\to X$ be the splitting covering of $f$. Then $E_f$ is connected.
Since $f$ splits on the identity covering $X\to X$, by the minimality property of the
splitting covering there is a covering map $\phi:X\to E_f$ over $X$.
But $X$ is connected and $\phi$ is a covering map onto the connected space $E_f$; hence
$\phi$ is surjective and has discrete fibres. Because $q\circ \phi=\mathrm{id}_X$, it follows that
$q$ is injective; therefore $q$ is a covering isomorphism and $E_f\cong X$.

Thus $\Deck(E_f/X)=1$ for every $f$, so the inverse limit $\PiST(X,x)$ is trivial.
Finally, continuous cohomology of the trivial group vanishes in positive degrees,
so $\HST^n(X,A)=0$ for $n>0$.
\end{proof}

\begin{corollary}\label{cor:torsion}
If $\pi_1(X, x)$ is a torsion group (not necessarily finite), then $\PiST(X, x)=1$ and
$\HST^n(X,A)=0$ for all $n>0$ and all discrete coefficient modules $A$.
\end{corollary}

\begin{proof}
A torsion group admits no nontrivial homomorphism to a torsion-free group, so it is
braid-free. Apply Theorem~\ref{thm:braidfree-trivial}.
\end{proof}

\begin{corollary}\label{cor:K(G,1)}
Let $G$ be a torsion group and let $X=K(G,1)$ be any Eilenberg--Mac Lane space of type $(G, 1)$.
Let $A$ be a finite abelian group with trivial action.
Then $\HST^2(X,A)=0$. In particular if $H^2(G,A)\neq 0$, then the comparison map $\Phi_X^2$ defined in Definition \ref{def:comparison},
\[
\Phi_X^2:\HST^2(X,A)\longrightarrow H^2_{\mathrm{sing}}(X,A)\cong H^2(G,A)
\]
is not surjective.
\end{corollary}

\begin{proof}
By Corollary~\ref{cor:torsion}, $\PiST(X, x)=1$, hence $\HST^2(X,A)=0$.
Since $X=K(G,1)$, group cohomology identifies with singular cohomology:
$H^2_{\mathrm{sing}}(X,A)\cong H^2(G,A)$. If $H^2(G,A)\neq 0$, a map from the zero group
cannot be surjective.
\end{proof}

\subsection{Cohomological dimension and vanishing}

\begin{example}
For $X=\RP^2$, we have $\pi_1(X, x)\cong \Z/2$, hence Corollary~\ref{cor:torsion} gives
$\PiST(\RP^2)=1$ and $\HST^n(\RP^2,A)=0$ for all $n>0$.
In contrast, $H^2_{\mathrm{sing}}(\RP^2,\Z/2)\cong \Z/2$, so the comparison map
$\Phi^2_X$ is the zero map and the unique nontrivial class is not Weierstrass realizable.
\end{example}

Recall that the \emph{cohomological dimension} of a profinite group
$G$, denoted by $\operatorname{cd}(G)$, is the smallest integer $d$
such that
\[
H^n_{\mathrm{cont}}(G,A)=0
\]
for every $n>d$ and every discrete torsion $G$-module $A$.
The \emph{cohomological dimension} of a topological space
$X$, denoted by $\operatorname{cd}(X)$, is the smallest integer $d$
such that
\[
H^n_{\mathrm{sing}}(X,A)=0
\]
for every $n>d$ and every coefficient group $A$.

\begin{definition}
The \emph{semi-topological cohomological dimension} of $X$ is
\[
\mathrm{cd}_{\mathrm{ST}}(X):=\mathrm{cd}(\PiST(X, x)),
\]
the smallest $d$ such that $\HST^n(X,A)=0$ for all $n>d$ and all discrete torsion
$\PiST(X, x)$-modules $A$.
\end{definition}

\begin{example}
Let $X=\mathbb{RP}^2$. Then its semi-topological cohomological dimension
$\mathrm{cd}_{ST}(X)$ is different from its (classical) cohomological dimension
$\mathrm{cd}(X)$:
\[
\mathrm{cd}_{ST}(\mathbb{RP}^2)=0
\qquad\text{but}\qquad
\mathrm{cd}(\mathbb{RP}^2)=2.
\]

More generally: any space $X$ with finite $\pi_1(X, x)$ but with some nonzero
$H^n_{\mathrm{sing}}(X;A)$ for $n\ge 1$ will force
$\mathrm{cd}_{ST}(X)=0$ by Corollary~\ref{cor:torsion}, while $\mathrm{cd}(X)\ge n$.
\end{example}

\subsection{Splitting coverings on the torus}

Let $X=T^2$ denote the two-dimensional torus. Then $\pi_1(X, x)\cong \Z^2$.
Any Weierstrass polynomial $f$ determines a braid monodromy
$\theta_f:\Z^2\to B_n$, i.e.\ a pair of commuting braids.
This imposes constraints on which finite quotients of $\Z^2$ occur as deck groups of
splitting coverings.

For simplicity of notation, we suppress the base point of $X$ in what follows.
\begin{theorem}\label{thm:torus-surj}
Let $X=T^2=S^1\times S^1\subset \C^2$ and let $A$ be a finite abelian group with trivial action.
The comparison map
\[
\Phi_X^2:\HST^2(T^2,A)\longrightarrow H^2_{\mathrm{sing}}(T^2,A)
\]
is surjective.
\end{theorem}

\begin{lemma}\label{lem:basic-cyclic}
Let
\[
u,v:T^2\to \C,\qquad u(e^{i\theta},e^{i\phi})=e^{i\theta},\quad v(e^{i\theta},e^{i\phi})=e^{i\phi}
\]
be the coordinate functions.
For $m\ge 1$, the polynomials
\[
f_1(z)=z^m-u,\qquad f_2(z)=z^m-v
\]
are Weierstrass on $T^2$, and their splitting coverings are the standard cyclic covers
\[
p_1:E_1=S^1\times S^1\to T^2,\quad (s,t)\mapsto (s^m,t),
\]
\[
p_2:E_2=S^1\times S^1\to T^2,\quad (s,t)\mapsto (s,t^m),
\]
with deck groups $\Z/m$ (acting by $s\mapsto \zeta s$ and $t\mapsto \zeta t$).
\end{lemma}

\begin{proof}
For each $(u,v)\in T^2$, the polynomial $z^m-u$ has $m$ distinct roots since $u\neq 0$,
so it is a Weierstrass polynomial. Its solution space
\(
\{((u,v),z):z^m=u\}
\)
identifies with $\{(s,t):|s|=|t|=1\}$ via $u=s^m$, giving the cover $p_1$.
On this space the pullback splits with roots $\zeta^k s$.

For minimality: if $q:Y\to T^2$ is any covering on which $z^m-u$ splits, choose one
continuous root $r\in C(Y)$ with $r^m=u\circ q$; then
\(y\mapsto (r(y),\mathrm{pr}_2(q(y)))\)
defines a covering map $Y\to E_1$ over $T^2$, showing $E_1$ is the splitting covering.
The case of $f_2$ is identical.
\end{proof}

\begin{lemma}\label{lem:dominating-product}
Fix any real constant $M>2$ and define
\[
h(z):=(z^m-u)\cdot((z-M)^m-v)\in C(T^2)[z].
\]
Then $h$ is a Weierstrass polynomial and its splitting covering $E_h\to T^2$ dominates both
$E_1\to T^2$ and $E_2\to T^2$.
\end{lemma}

\begin{proof}
For any $(u,v)\in T^2$, all roots of $z^m-u$ lie on the unit circle, hence in the closed disk
$\overline{D(0,1)}$. All roots of $(z-M)^m-v$ are of the form $M+\eta$ with $|\eta|=1$,
hence lie in $\overline{D(M,1)}$. If $M>2$, these disks are disjoint, so $h$ has $2m$
distinct roots and is a Weierstrass polynomial.

Let $p:E_h\to T^2$ be the splitting covering of $h$. In $E_h$ the polynomial $p^*h$ splits
as a product of linear factors with continuous roots $r_1,\dots,r_{2m}\in C(E_h)$.
Because $\overline{D(0,1)}$ and $\overline{D(M,1)}$ are disjoint, the subset of indices
$i$ for which $r_i(\cdot)\in\overline{D(0,1)}$ is locally constant, hence constant on the
connected space $E_h$. Exactly $m$ roots lie in $\overline{D(0,1)}$ and their product is a
monic degree--$m$ factor of $p^*h$ whose fibrewise roots are precisely the roots of
$p^*(z^m-u)$; thus $z^m-u$ splits in $E_h$. By Lemma~\ref{lem:basic-cyclic} and
minimality of the splitting covering, there is a covering map $E_h\to E_1$ over $T^2$.
The same argument with $\overline{D(M,1)}$ shows $z^m-v$ splits in $E_h$, giving
$E_h\to E_2$.
\end{proof}

\begin{corollary}\label{cor:ZmxZm-quot}
The multiplication-by-$m$ covering
\[
[m]:T^2\to T^2,\qquad (s,t)\mapsto (s^m,t^m),
\]
whose deck group is $Q=(\Z/m)^2$, is dominated by the splitting covering $E_h\to T^2$.
Equivalently, there is a continuous surjection
\[
\PiST(T^2)\twoheadrightarrow Q=(\Z/m)^2.
\]
\end{corollary}

\begin{proof}
The fibre product $E_1\times_{T^2}E_2\to T^2$ is precisely $[m]:T^2\to T^2$ and has deck
group $(\Z/m)^2$. Since $E_h$ dominates both $E_1$ and $E_2$ (Lemma~\ref{lem:dominating-product}),
it dominates their fibre product, hence dominates $[m]$. By definition of $\PiST(T^2)$
as an inverse limit over deck groups of splitting coverings, this means $(\Z/m)^2$ occurs
as a quotient of $\PiST(T^2)$.
\end{proof}

\begin{proof}[Proof of Theorem 4.13]
Since $T^2\simeq K(\Z^2,1)$,
$
H^2_{\mathrm{sing}}(T^2,A)\cong H^2(\Z^2,A).
$
Also $H_2(T^2;\Z)\cong \Z$, so the universal coefficient theorem gives
\[
H^2_{\mathrm{sing}}(T^2,A)\cong \Hom(H_2(T^2;\Z),A)\cong \Hom(\Z,A)\cong A.
\]
Thus it suffices to show: for each $a\in A$, the corresponding class in $H^2_{\mathrm{sing}}(T^2,A)$
lies in the image of $\Phi_X^2$.

\medskip
\textbf{Step 1: a canonical class factoring through $(\Z/m)^2$.}
Let $m:=\exp(A):=\min\{n\geq 1: na=0 \mbox{ for all } a\in A\}$ be the exponent of $A$ and fix a homomorphism $\lambda:\Z/m\to A$ with $\lambda(1)=a$.
Let $Q:=(\Z/m)^2$ and define a normalized $2$--cochain $c:Q\times Q\to A$ by
\[
c\big((x_1,x_2),(y_1,y_2)\big):=\lambda(x_1y_2)
\]
where arithmetic is in $\Z/m$.
This is a standard cocycle representing the alternating bilinear form (``area form'')
on $Q$, hence determines a class $[c]\in H^2(Q,A)$.

Let $\pi:\Z^2\twoheadrightarrow Q$ be reduction mod $m$ in each coordinate. Then the
pullback $\pi^*[c]\in H^2(\Z^2,A)\cong H^2_{\mathrm{sing}}(T^2,A)$ is precisely the class corresponding
to $a\in A$ under $H^2_{\mathrm{sing}}(T^2,A)\cong A$ (equivalently, it is $a$ times the fundamental class).

Therefore, to realize the class corresponding to $a$, it is enough to prove that the
quotient $\pi:\pi_1(T^2)=\Z^2\twoheadrightarrow Q$ occurs as a quotient of $\PiST(T^2)$,
i.e.\ that the covering $[m]:T^2\to T^2$ with deck group $Q=(\Z/m)^2$ is dominated by
a splitting covering arising from a Weierstrass polynomial on $T^2$.

\medskip
\textbf{Step 2: build a class in $\HST^2$ mapping to the chosen $a\in A$.}
Let $q:\PiST(T^2)\twoheadrightarrow Q$ be the surjection from Corollary~\ref{cor:ZmxZm-quot}.
Inflate the class $[c]\in H^2(Q,A)$ along $q$:
\[
\beta:=q^*[c]\in H^2_{\cont}(\PiST(T^2),A)=\HST^2(T^2,A).
\]
By construction, the image $\Phi_X^2(\beta)\in H^2(T^2,A)\cong H^2(\Z^2,A)$ is the pullback
of $[c]$ along the composite $\Z^2\to \PiST(T^2)\to Q$, which coincides with
$\pi:\Z^2\twoheadrightarrow Q$ (it is the monodromy of the cover $[m]$).
Hence
\[
\Phi_X^2(\beta)=\pi^*[c],
\]
and by Step 1 this is exactly the class corresponding to $a\in A$.

Since $a\in A$ was arbitrary, every class in $H^2_{\mathrm{sing}}(T^2,A)$ lies in the image of $\Phi_X^2$.
Therefore $\Phi_X^2$ is surjective.
\end{proof}

\subsection{Linearizing finite projective monodromy}

We now apply semi-topological Galois cohomology to the problem of
lifting finite projective monodromy to genuine linear monodromy.

\begin{definition}
A \emph{finite projective semi-topological local system of rank $n$}
on $X$ is a continuous homomorphism
\[
\overline{\rho}:
\Pi^{\mathrm{ST}}(X, x)
\longrightarrow
\operatorname{PGL}_n(\mathbb{C})
\]
whose image is finite.

We say that $\overline{\rho}$ is
\emph{semi-topologically linearizable} if there exists a continuous
homomorphism
\[
\rho:
\Pi^{\mathrm{ST}}(X, x)
\longrightarrow
\operatorname{GL}_n(\mathbb{C})
\]
such that
\[
P\circ\rho=\overline{\rho},
\]
where
\[
P:
\operatorname{GL}_n(\mathbb{C})
\twoheadrightarrow
\operatorname{PGL}_n(\mathbb{C})
\]
is the canonical projectivization map. Thus semi-topological
linearizability means that the projective monodromy is induced by an
actual linear representation of $\Pi^{\mathrm{ST}}(X, x)$.
\end{definition}

Let
\[
G:=\operatorname{im}(\overline{\rho})
\subset \operatorname{PGL}_n(\mathbb{C}).
\]
We recall the following standard fact from the theory of projective
representations of finite groups; see \cite{SerreLRFG}.

\begin{lemma}
Let $G\subset \operatorname{PGL}_n(\mathbb{C})$
be a finite subgroup. Then there exist an integer $m\geq 1$ and a
finite subgroup
$\widetilde{G}\subset \operatorname{GL}_n(\mathbb{C})$
such that the projectivization map restricts to a surjective
homomorphism
\[
\pi:=P|_{\widetilde{G}}:
\widetilde{G}\twoheadrightarrow G
\]
and fits into a central extension
\[
1\longrightarrow \mu_m
\longrightarrow \widetilde{G}
\stackrel{\pi}{\longrightarrow}G
\longrightarrow1,
\]
where
\[
\mu_m=
\{\zeta\in\mathbb{C}^{\times}:\zeta^m=1\}
\]
is embedded as scalar matrices. The extension determines a class
\[
\alpha_{\widetilde{G}}
\in H^2(G,\mu_m),
\]
where $G$ acts trivially on $\mu_m$.
\end{lemma}

The finite central lift $\widetilde{G}$ need not be unique, and hence
the class $\alpha_{\widetilde{G}}$ may depend on the chosen lift.

Let
\[
q:
\Pi^{\mathrm{ST}}(X, x)
\twoheadrightarrow G
\]
be the quotient homomorphism induced by $\overline{\rho}$. Since
continuous group cohomology is contravariant, $q$ induces the
inflation map
\[
q^*:
H^2(G,\mu_m)
\longrightarrow
H^2_{\mathrm{cont}}
\bigl(\Pi^{\mathrm{ST}}(X, x),\mu_m\bigr)
=
H^2_{\mathrm{ST}}(X,\mu_m).
\]
Thus a finite central lift $\widetilde{G}$ determines an obstruction
class
\[
q^*(\alpha_{\widetilde{G}})
\in H^2_{\mathrm{ST}}(X,\mu_m).
\]

We refer the reader to Serre's book \cite[Chapter I, Proposition 43]{SerreGaloisCohom} for a result related to the following.

\begin{theorem}
Let
\[
\overline{\rho}:
\Pi^{\mathrm{ST}}(X, x)
\longrightarrow
\operatorname{PGL}_n(\mathbb{C})
\]
be a finite projective semi-topological local system, and let
$G:=\operatorname{im}(\overline{\rho})$.
Then the following are equivalent:

\begin{enumerate}
\item[(i)]
The projective semi-topological local system $\overline{\rho}$ is
semi-topologically linearizable.

\item[(ii)]
There exist an integer $m\geq1$ and a finite subgroup
$\widetilde{G}\subset\operatorname{GL}_n(\mathbb{C})$
such that
\[
1\longrightarrow\mu_m
\longrightarrow\widetilde{G}
\stackrel{\pi}{\longrightarrow}G
\longrightarrow1
\]
is a central extension, where
$\pi=P|_{\widetilde{G}}$,
and the quotient homomorphism
$q:\Pi^{\mathrm{ST}}(X, x)\twoheadrightarrow G$
admits a continuous lift
$\widetilde{\rho}:\Pi^{\mathrm{ST}}(X, x)\longrightarrow\widetilde{G}$
satisfying
\[
\pi\circ\widetilde{\rho}=q
\]

\item[(iii)]
There exist an integer $m\geq1$ and a finite subgroup
$\widetilde{G}\subset\operatorname{GL}_n(\mathbb{C})$
such that
\[
1\longrightarrow\mu_m
\longrightarrow\widetilde{G}
\stackrel{\pi}{\longrightarrow}G
\longrightarrow1,
\qquad
\pi=P|_{\widetilde{G}},
\]
is a central extension and
\[
q^*(\alpha_{\widetilde{G}})=0
\qquad\text{in}\qquad
H^2_{\mathrm{ST}}(X,\mu_m).
\]
\end{enumerate}
\end{theorem}

\begin{proof}
We first prove $(i)\Rightarrow(ii)$. Suppose that
$\overline{\rho}$ is semi-topologically linearizable, and choose a
continuous lift
\[
\rho:
\Pi^{\mathrm{ST}}(X, x)
\longrightarrow
\operatorname{GL}_n(\mathbb{C})
\]
such that
\[
P\circ\rho=\overline{\rho}.
\]
Set
\[
\widetilde{G}:=\operatorname{im}(\rho).
\]

Since $\Pi^{\mathrm{ST}}(X, x)$ is profinite and $\rho$ is continuous,
$\widetilde{G}$ is a compact totally disconnected subgroup of
$\operatorname{GL}_n(\mathbb{C})$. It is therefore a closed subgroup
of the Lie group $\operatorname{GL}_n(\mathbb{C})$. By the closed
subgroup theorem, $\widetilde{G}$ is itself a Lie group. Since it is
totally disconnected, it has dimension zero and is therefore
discrete. Being both compact and discrete, $\widetilde{G}$ is finite.

The projectivization map restricts to a surjection
\[
P|_{\widetilde{G}}:
\widetilde{G}\twoheadrightarrow G.
\]
Its kernel consists precisely of those elements of $\widetilde{G}$
which are scalar matrices. Hence
$\ker(P|_{\widetilde{G}})$
is a finite subgroup of $\mathbb{C}^{\times}$ and is therefore cyclic.
Thus, for some $m\geq1$,
\[
\ker(P|_{\widetilde{G}})\cong\mu_m.
\]
We obtain a finite central extension
\[
1\longrightarrow\mu_m
\longrightarrow\widetilde{G}
\stackrel{\pi}{\longrightarrow}G
\longrightarrow1,
\]
where $\pi=P|_{\widetilde{G}}$. The homomorphism $\rho$, regarded as
a map into $\widetilde{G}$, gives the required lift of $q$. Hence
$(i)$ implies $(ii)$.

Conversely, assume (ii). Composing
$\widetilde{\rho}:\Pi^{\mathrm{ST}}(X, x)\longrightarrow\widetilde{G}$
with the inclusion
$\widetilde{G}\hookrightarrow\operatorname{GL}_n(\mathbb{C})$
gives a continuous homomorphism
\[
\rho:
\Pi^{\mathrm{ST}}(X, x)\longrightarrow
\operatorname{GL}_n(\mathbb{C}).
\]
Since $\pi=P|_{\widetilde{G}}$ and
$\pi\circ\widetilde{\rho}=q$, its projectivization is
$\overline{\rho}$. Hence $(ii)$ implies $(i)$.

It remains to prove $(ii)\Longleftrightarrow(iii)$. Fix a finite
central extension
\[
1\longrightarrow\mu_m
\longrightarrow\widetilde{G}
\stackrel{\pi}{\longrightarrow}G
\longrightarrow1
\]
as in (ii), and let
$\alpha_{\widetilde{G}}\in H^2(G,\mu_m)$
be its extension class.

Pulling this extension back along
$q:\Pi^{\mathrm{ST}}(X, x)\twoheadrightarrow G$
gives the central extension
\[
1\longrightarrow\mu_m
\longrightarrow
\Pi^{\mathrm{ST}}(X, x)\times_G\widetilde{G}
\longrightarrow
\Pi^{\mathrm{ST}}(X, x)
\longrightarrow1,
\]
where
\[
\Pi^{\mathrm{ST}}(X, x)\times_G\widetilde{G}
=
\bigl\{
(\gamma,\widetilde{g})\mid
q(\gamma)=\pi(\widetilde{g})
\bigr\}
\]
The cohomology class of this pulled-back extension is precisely
\[
q^*(\alpha_{\widetilde{G}})
\in
H^2_{\mathrm{cont}}
\bigl(\Pi^{\mathrm{ST}}(X, x),\mu_m\bigr)
=
H^2_{\mathrm{ST}}(X,\mu_m).
\]

By the standard classification of continuous central extensions by
second continuous group cohomology, this class vanishes if and only
if the pulled-back extension splits; see \cite{SerreGaloisCohom}. A splitting is
equivalent to a continuous homomorphism
\[
\widetilde{\rho}:
\Pi^{\mathrm{ST}}(X, x)
\longrightarrow\widetilde{G}
\]
such that
\[
\pi\circ\widetilde{\rho}=q.
\]
Therefore $(ii)$ and $(iii)$ are equivalent.
\end{proof}

\begin{remark}[Interpretation]
A projective representation assigns to each monodromy element a linear
transformation only up to multiplication by a nonzero scalar. Thus a
$\operatorname{PGL}_n(\mathbb{C})$-local system records projective
linear monodromy, whereas a \emph{genuine}
$\operatorname{GL}_n(\mathbb{C})$-local system has monodromy given by
an actual linear representation
$\rho:\Pi^{\mathrm{ST}}(X, x)\longrightarrow\operatorname{GL}_n(\mathbb{C})$.

To describe the obstruction more concretely, choose matrix
representatives
\[
A_g\in\operatorname{GL}_n(\mathbb{C}),
\qquad g\in G,
\]
for the projective elements of $G$. In general these representatives
do not satisfy the group law exactly. Instead,
\[
A_gA_h=c(g,h)A_{gh},
\]
where
$c(g,h)\in\mathbb{C}^{\times}$ and $\C^{\times}$ is the multiplicative group of nonzero complex numbers.
Associativity implies that $c$ is a $2$-cocycle, where
$\mathbb{C}^{\times}$ is regarded as a trivial $G$-module. Its
cohomology class
\[
[c]\in H^2(G,\mathbb{C}^{\times})
\]
is called the \emph{factor-set class} of the projective
representation. Its vanishing is equivalent to the possibility of
choosing the representatives $A_g$ so that they form an ordinary
linear representation of $G$.

Recall that the \emph{Schur multiplier} of a finite group $G$ is
\[
M(G):=H_2(G,\mathbb{Z}).
\]
Since $\mathbb{C}^{\times}$ is a divisible abelian group, the universal
coefficient theorem gives a natural isomorphism
\[
H^2(G,\mathbb{C}^{\times})
\cong
\operatorname{Hom}
\bigl(M(G),\mathbb{C}^{\times}\bigr).
\]
Thus the factor-set obstruction of a projective representation is
naturally controlled by the Schur multiplier.

For finite $G$, the factor-set class has finite order. Consequently,
for some $m\geq1$, it may be represented by a finite central extension
\[
1\longrightarrow\mu_m
\longrightarrow\widetilde{G}
\longrightarrow G
\longrightarrow1
\]
with
$\widetilde{G}\subset\operatorname{GL}_n(\mathbb{C})$.
The associated class
$\alpha_{\widetilde{G}}\in H^2(G,\mu_m)$
pulls back along
$q:\Pi^{\mathrm{ST}}(X, x)\twoheadrightarrow G$
to
\[
q^*(\alpha_{\widetilde{G}})
\in H^2_{\mathrm{ST}}(X,\mu_m).
\]

Theorem~4.19 shows that the vanishing of this pulled-back class is
precisely the obstruction to lifting the finite projective
semi-topological monodromy through the chosen central extension.
Allowing the choice of a finite central linear lift
$\widetilde{G}$, the projective semi-topological local system
$\overline{\rho}$ is semi-topologically linearizable if and only if
such an obstruction class vanishes.
\end{remark}


\section{The $\pi_1$-detectable Weierstrass realizability conjecture}

Let $X$ be a smooth complex projective variety and fix an integer $m\ge 2$.
Write $A:=\mathbb Z/m$ with trivial action.

Recall that the \emph{N\'eron--Severi group} of $X$ is
\[
\operatorname{NS}(X)
:=
\operatorname{Pic}(X)/\operatorname{Pic}^0(X),
\]
where $\operatorname{Pic}(X)$ is the Picard group of line bundles on
$X$ and $\operatorname{Pic}^0(X)$ is the subgroup of line bundles
algebraically equivalent to zero.

The first Chern class induces a homomorphism
\[
c_1:
\operatorname{NS}(X)
\longrightarrow
H^2_{\mathrm{sing}}(X,\mathbb{Z}).
\]
For $m\geq 2$, define
\[
\rho_m:
\operatorname{NS}(X)
\xrightarrow{\;c_1\;}
H^2_{\mathrm{sing}}(X,\mathbb{Z})
\longrightarrow
H^2_{\mathrm{sing}}(X,\mathbb{Z}/m),
\]
where the second homomorphism is induced by the coefficient reduction
$\mathbb{Z}\longrightarrow\mathbb{Z}/m$.
Thus $\rho_m$ sends the class of a divisor, or equivalently the class
of a line bundle, to the reduction modulo $m$ of its first Chern
class.

Recall that in the semi-topological setting one has a canonical degree-$2$ comparison map defined in Definition \ref{def:comparison},
\[
\Phi_X^2:\ H^2_{ST}(X,\mathbb Z/m)\longrightarrow H^2_{\mathrm{sing}}(X,\mathbb Z/m),
\]
as the composite
\[
H^2_{ST}(X,A)=H^2_{\mathrm{cont}}(\PiST(X, x),A)
\xrightarrow{\ \rho^*\ } H^2_{\mathrm{cont}}(\widehat{\pi_1(X, x)},A)
\xrightarrow{\ \iota^*\ } H^2(\pi_1(X, x),A)
\xrightarrow{\ c^*\ } H^2_{\mathrm{sing}}(X,A),
\]
and a class $\alpha\in H^2_{\mathrm{sing}}(X,A)$ is \emph{Weierstrass realizable} if
$\alpha\in \mathrm{Im}(\Phi_X^2)$.

\subsection{$\pi_1$-detectable classes}

Define the $\pi_1$-detectable subgroup
\[
H^2_{\mathrm{sing}}(X,A)_{\pi_1}
:=\mathrm{Im}\!\Bigl(c^*:H^2(\pi_1(X, x),A)\to H^2_{\mathrm{sing}}(X,A)\Bigr).
\]
By construction of $\Phi_X^2$ one always has the necessary containment
\[
\mathrm{Im}(\Phi_X^2)\subseteq H^2_{\mathrm{sing}}(X,A)_{\pi_1}.
\]
Hence only $\pi_1$-detectable degree-$2$ classes can possibly be Weierstrass realizable.

\begin{conjecture}[The $\pi_1$-detectable Weierstrass realizability conjecture]\label{conj:divisor-pi1}
For every smooth complex projective variety $X$ and every $m\ge 2$,
\[
\rho_m(NS(X))\ \cap\ H^2_{\mathrm{sing}}(X,\mathbb Z/m)_{\pi_1}
=\rho_m(NS(X))\ \cap \mathrm{Im}(\Phi_X^2)
\]
Equivalently: every divisor class modulo $m$ that is $\pi_1$-detectable is Weierstrass realizable.
\end{conjecture}

\begin{remark}
If $X$ is aspherical (a $K(\pi,1)$), then $c^*$ is an isomorphism in all degrees. In that case
\[
H^2_{\mathrm{sing}}(X,\mathbb Z/m)_{\pi_1} = H^2_{\mathrm{sing}}(X,\mathbb Z/m),
\]
and Conjecture~\ref{conj:divisor-pi1} reduces to the simpler statement
\[
\rho_m(NS(X))\subseteq \mathrm{Im}(\Phi_X^2).
\]
\end{remark}

If $H^2(\pi_1(X, x),\mathbb Z/m)=0$, then $H^2_{\mathrm{sing}}(X,\mathbb Z/m)_{\pi_1}=0$ and the
left-hand side of Conjecture~\ref{conj:divisor-pi1} is zero. In particular, if $X$ is simply
connected, the conjecture holds vacuously.

\subsection{Abelian varieties}

\begin{theorem}\label{thm:abelian-varieties}
Let $X$ be a complex abelian variety of dimension $g$ and fix $m\ge 2$.
Put $A:=\mathbb Z/m$ with trivial action. Then the degree-$2$ comparison map
\[
\Phi_X^2:\ H^2_{ST}(X,A)\longrightarrow H^2_{\mathrm{sing}}(X,A)
\]
is surjective. In particular, the $\pi_1$-detectable Weierstrass realizability conjecture is true for abelian varieties.
\end{theorem}

\begin{proof}
\textbf{Step 1: Reduce to the real torus.}
As a complex torus, $X$ is of the form $X\simeq \mathbb C^g/\Lambda$ for a lattice
$\Lambda\simeq \mathbb Z^{2g}$. Choosing an $\mathbb R$-basis identifies
$\mathbb C^g\simeq \mathbb R^{2g}$ and $\Lambda\simeq \mathbb Z^{2g}$, hence yields a
homeomorphism
\[
h:\ T^{2g}=\mathbb R^{2g}/\mathbb Z^{2g}\ \xrightarrow{\ \sim\ }\ X.
\]
By contravariant functoriality of $H^\bullet_{ST}$, the pullback $h^*$
induces isomorphisms
\[
h^*:H^2_{ST}(X,A)\xrightarrow{\sim} H^2_{ST}(T^{2g},A),
\qquad
h^*:H^2_{\mathrm{sing}}(X,A)\xrightarrow{\sim} H^2_{\mathrm{sing}}(T^{2g},A),
\]
and by Theorem \ref{thm:comparison-naturality}, the functoriality of the construction of $\Phi^2$ shows that
\[
h^*\circ \Phi_X^2 \;=\; \Phi_{T^{2g}}^2\circ h^*.
\]
Therefore it suffices to prove that $\Phi_{T^{2g}}^2$ is surjective.

\medskip
\textbf{Step 2: A generating set for $H^2_{\mathrm{sing}}(T^N,A)$.}
Let $N:=2g$. For $1\le i<j\le N$, let
\[
p_{ij}:T^N\to T^2
\]
be the projection onto the $(i,j)$-coordinates.
Let $\omega\in H^2_{\mathrm{sing}}(T^2,A)$ denote the class corresponding to $1\in A$
under the standard identification $H^2_{\mathrm{sing}}(T^2,A)\cong A$.
Then the classes $p_{ij}^*(a\omega)$, as $i<j$ vary and $a\in A$ varies, generate
$H^2_{\mathrm{sing}}(T^N,A)$.

Indeed, $H^\bullet_{\mathrm{sing}}(T^N,\mathbb Z)\cong \bigwedge^\bullet \mathbb Z^N$ is torsion-free, so by the
universal coefficient theorem
\[
H^2_{\mathrm{sing}}(T^N,A)\ \cong\ H^2_{\mathrm{sing}}(T^N,\mathbb Z)\otimes A \ \cong\ \bigwedge\nolimits^2\mathbb Z^N\otimes A,
\]
and the wedge basis elements $e_i\wedge e_j$ correspond exactly to the pullbacks
$p_{ij}^*(\omega)$.

\medskip
\textbf{Step 3: Surjectivity of $\Phi_{T^N}^2$ by reduction to $T^2$.}
Fix $i<j$ and $a\in A$.
By Theorem \ref{thm:torus-surj}, there exists $\xi_{a}\in H^2_{ST}(T^2,A)$ such that
\[
\Phi_{T^2}^2(\xi_a)=a\omega.
\]
Pulling back along $p_{ij}$ gives an element $p_{ij}^*(\xi_a)\in H^2_{ST}(T^N,A)$. 
By Theorem \ref{thm:comparison-naturality}, the naturality of $\Phi^2$ with respect to pullback gives
\[
\Phi_{T^N}^2\bigl(p_{ij}^*(\xi_a)\bigr)
\;=\;
p_{ij}^*\bigl(\Phi_{T^2}^2(\xi_a)\bigr)
\;=\;
p_{ij}^*(a\omega).
\]
Thus each generator $p_{ij}^*(a\omega)$ lies in $\mathrm{Im}(\Phi_{T^N}^2)$, hence
$\Phi_{T^N}^2$ is surjective.

\medskip
\textbf{Step 4: Conclude for $X$ and divisor classes.}
Surjectivity for $T^{2g}$ implies surjectivity for $X$ by Step~1. Hence
$\rho_m(NS(X))\subseteq H^2_{\mathrm{sing}}(X,A)=\mathrm{Im}(\Phi_X^2)$, proving the theorem.
\end{proof}

\subsection{Smooth complex projective curves}
\begin{theorem}
Let $C$ be a smooth complex projective curve of genus $g\ge 1$, and let $m\ge 2$.
Put $A:=\Z/m$ with trivial action. Then the degree-$2$ comparison map
\[
\Phi_C^2:H^2_{ST}(C,A)\longrightarrow H^2_{\mathrm{sing}}(C,A)
\]
is surjective. Equivalently, every class in $H^2_{\mathrm{sing}}(C,\Z/m)$ is Weierstrass realizable.
In particular, the $\pi_1$-detectable Weierstrass realizability conjecture is true for $C$.
\end{theorem}

\begin{remark}
The statement fails for $C=\PP^1$. Indeed, $\pi_1(\PP^1)=1$, so $H^2_{ST}(\PP^1,A)=0$
whereas $H^2_{\mathrm{sing}}(\PP^1,A)\cong A\neq 0$.
Thus the surjectivity is a genus $\ge 1$ phenomenon. Note that the $\pi_1$-detectable Weierstrass realizability conjecture is trivially
true for $\CP^n$ for any $n\in \N$.
\end{remark}

\begin{proof}
Since $C$ is a closed orientable surface of genus $g$, it is a $K(\pi_1(C, c_0),1)$. In particular,
\[
H^2_{\mathrm{sing}}(C,A)\cong H^2(\pi_1(C, c_0),A)\cong A,
\]
so it suffices to produce one class in $H^2_{ST}(C,A)$ mapping to a generator of
$H^2_{\mathrm{sing}}(C,A)$.

\medskip
\noindent\textbf{Step 1: choose $2g$ circle-valued functions realizing the abelianization.}
Since
\[
H^1_{\mathrm{sing}}(C,\Z)\cong \mathrm{Hom}(\pi_1(C, c_0),\Z)\cong \Z^{2g},
\]
and $(S^1)^{2g}$ is a $K(\Z^{2g},1)$, there exists a continuous map
$\alpha:C\longrightarrow (S^1)^{2g}$
inducing the abelianization map
\[
\pi_1(C, c_0) \twoheadrightarrow H_1(C,\Z)\cong \Z^{2g}.
\]
Write the coordinate functions of $\alpha$ as
$
u_1,\dots,u_{2g}:C\to S^1\subset \C.
$
Then the cohomology classes $[u_1],\dots,[u_{2g}]\in H^1_{\mathrm{sing}}(C,\Z)$ form a basis.
Reducing modulo $m$, they induce a surjective homomorphism
$\chi:\pi_1(C, c_0)\twoheadrightarrow (\Z/m)^{2g}$.
Let $Q:=(\Z/m)^{2g}$.

\medskip
\noindent\textbf{Step 2: the cyclic covers defined by $z^m-u_i$ are semi-topological.}

For each $i=1,\dots,2g$, consider the polynomial
\[
f_i(z):=z^m-u_i\in C(C)[z].
\]
Since $u_i(x)\in S^1$ for every $x$, the polynomial $z^m-u_i(x)$ has $m$ distinct roots, so
$f_i$ is a Weierstrass polynomial.

Let $\chi_i:\pi_1(C, c_0)\to \Z/m$ be the mod-$m$ reduction of the class $[u_i]\in H^1_{\mathrm{sing}}(C,\Z)$, and let
$p_i:E_i\to C$
be the connected covering corresponding to $\ker(\chi_i)$.
We claim that $E_i\to C$ is the splitting covering of $f_i$.
Indeed, since the $m$-th power map
\[
q_m:S^1\to S^1,\qquad z\mapsto z^m
\]
corresponds to the subgroup $m\Z\subset \Z=\pi_1(S^1)$, the pullback $u_i\circ p_i:E_i\to S^1$
lifts to a map $r_i:E_i\to S^1$
such that
\[
r_i^m=u_i\circ p_i.
\]
Hence $f_i$ splits on $E_i$:
\[
p_i^*f_i(z)=z^m-r_i^m=\prod_{k=0}^{m-1}(z-\zeta^k r_i),
\qquad \zeta=e^{2\pi i/m}.
\]

Conversely, if $Y\to C$ is any connected covering on which $f_i$ splits, then there exists a
continuous root
\[
s:Y\to S^1,\qquad s^m=u_i\circ p.
\]
Therefore $u_i\circ p$ lifts through $q_m$, so $p_*\pi_1(Y)\subseteq \ker(\chi_i)$.
By the covering-space classification, $Y\to C$ factors through $E_i\to C$.
Thus $E_i\to C$ is indeed the splitting covering of $f_i$.

\medskip
\noindent\textbf{Step 3: build one splitting covering dominating the full $Q$-cover.}

Choose real numbers $M_1,\dots,M_{2g}$
such that $|M_i-M_j|>2$ for $i\neq j$.
Define
\[
h(z):=\prod_{i=1}^{2g}\bigl((z-M_i)^m-u_i\bigr)\in C(C)[z].
\]
For each $x\in C$, the roots of $(z-M_i)^m-u_i(x)$ lie in the closed disk
$\overline{D(M_i,1)}$, and these disks are pairwise disjoint. Hence $h_x(z)$ has $2gm$
distinct roots for every $x$, so $h$ is a Weierstrass polynomial.

Let $p_h:E_h\to C$ be the splitting covering of $h$. Because the root clusters lie in pairwise
disjoint disks, on the connected space $E_h$ the $m$ roots lying in each disk \(\overline{D(M_i,1)}\)
vary continuously and define a splitting of the factor \((z-M_i)^m-u_i\). Thus each $f_i$ splits on $E_h$. By minimality of the
splitting covering $E_i\to C$, there are covering maps
$E_h\longrightarrow E_i$ for $i=1,\dots,2g$ over $C$.

Hence $E_h$ dominates the fibre product
\[
Y:=E_1\times_C\cdots\times_C E_{2g}.
\]
The covering $Y\to C$ is the connected cover corresponding to the subgroup
\[
\bigcap_{i=1}^{2g}\ker(\chi_i)=\ker(\chi),
\]
so its deck group is exactly
\[
\Deck(Y/C)\cong Q=(\Z/m)^{2g}.
\]
Therefore $Q$ occurs as a finite quotient of $\PiST(C)$. Equivalently, there is a canonical
continuous surjection
$\PiST(C)\twoheadrightarrow Q$.

\medskip
\noindent\textbf{Step 4: a degree-2 class on $Q$ giving the fundamental class of $C$.}

Let $x_1,y_1,\dots,x_g,y_g\in H^1(Q,A)$
be the coordinate classes corresponding to the standard basis of $Q=(\Z/m)^{2g}$. Define
\[
\omega_Q:=x_1\smile y_1 \in H^2(Q,A).
\]

Under the pullback
\[
\chi^*:H^2(Q,A)\to H^2(\pi,A)\cong H^2_{\mathrm{sing}}(C,A),
\]
the class $\omega_Q$ maps to the mod-$m$ fundamental class of $C$, hence to a generator of
$H^2_{\mathrm{sing}}(C,A)\cong A$.

Indeed, if $a_1,b_1,\dots,a_g,b_g\in H^1(C,A)$
is the basis dual to a symplectic basis of $H_1(C,\Z)$, then the mod-$m$ orientation class is
\[
[C]_m= a_1\smile b_1
\]

\medskip
\noindent\textbf{Step 5: inflate $\omega_Q$ to $H^2_{ST}(C,A)$.}

Let $q:\PiST(C)\twoheadrightarrow Q$
be the quotient map obtained in Step 3, and let
\[
\widetilde{\omega}:=q^*(\omega_Q)\in H^2_{\mathrm{cont}}(\PiST(C),A)=H^2_{ST}(C,A).
\]
By construction of the comparison map, $\Phi_C^2(\widetilde{\omega})$ is the pullback of
$\omega_Q$ along
\[
\pi_1(C, c_0)\longrightarrow \PiST(C)\xrightarrow{q}Q,
\]
which is exactly the homomorphism $\chi$. Therefore
\[
\Phi_C^2(\widetilde{\omega})=\chi^*(\omega_Q)=[C]_m.
\]
So the image of $\Phi_C^2$ contains a generator of $H^2_{\mathrm{sing}}(C,A)\cong A$.
Since the target is cyclic, $\Phi_C^2$ is surjective.
\end{proof}

\subsection{Ruled surfaces over positive-genus curves}
\begin{theorem}\label{thm:ruled-surface-pi1-detectable}
Let $C$ be a smooth complex projective curve of genus $g\ge 1$, let
\[
p:X=\PP(E)\longrightarrow C
\]
be a geometrically ruled surface over $C$, and let $m\ge 2$.
Write $A:=\Z/m$ with trivial action.

Then:
\begin{enumerate}[label=\textup{(\arabic*)}]
\item $X$ is not an Eilenberg--Mac Lane space;
\item the classifying map
\[
c_X^*:H^2(\pi_1(X, x),A)\longrightarrow H^2_{\mathrm{sing}}(X,A)
\]
is not an isomorphism;
\item one has
\[
\mathrm{Im}(\Phi_X^2)=H^2_{\mathrm{sing}}(X,A)_{\pi_1}=p^*H^2_{\mathrm{sing}}(C,A);
\]

\item The $\pi_1$-detectable Weierstrass realizability conjecture holds for every ruled surface
$\PP(E)\to C$ over a smooth projective curve of genus $g\ge 1$.
\end{enumerate}
\end{theorem}

\begin{proof}
Let $c_0$ be the base point of $C$.
Since $p:X\to C$ is a $\PP^1$-bundle and $\PP^1$ is simply connected, the long exact
sequence of homotopy groups yields an isomorphism
$p_*:\pi_1(X, x)\xrightarrow{\ \cong\ }\pi_1(C, c_0)$.
Moreover, the fibre $F\simeq \PP^1$ gives a nontrivial element of $\pi_2(X)$, so $X$ is not an Eilenberg--Mac Lane space proving \textup{(1)}.

Let $b:C\to B\pi_1(C, c_0)$ be a classifying map.
Since $C$ is a smooth complex projective curve of genus $g\geq 1$,
its underlying topological space is a closed 2-dimensional orientable surface. By the classification of orientable closed real surfaces,
\[
\pi_k(C)=0
\qquad\text{for all } k\geq 2
\]
Therefore $C$ is aspherical and has the homotopy type of an
Eilenberg--Mac Lane space $K(\pi_1(C, c_0),1)$.
Consequently, the classifying map
$b:C\longrightarrow B\pi_1(C, c_0)$ is a homotopy equivalence.
Via the isomorphism
$\pi_1(X, x)\cong \pi_1(C, c_0)$, we may choose the classifying map of $X$ to be
\[
c_X=b\circ p:X\longrightarrow B\pi_1(X, x)\cong B\pi_1(C, c_0).
\]
It follows that
\[
H^2_{\mathrm{sing}}(X,A)_{\pi_1}
=
\mathrm{Im}(c_X^*)
=
\mathrm{Im}(p^*\circ b^*)
=
p^*H^2_{\mathrm{sing}}(C,A).
\]
Thus
\[
H^2_{\mathrm{sing}}(X,A)_{\pi_1}=p^*H^2_{\mathrm{sing}}(C,A).
\tag{$*$}
\]

We next show that $c_X^*$ is not an isomorphism. Let
\[
\xi:=c_1(\mathcal O_X(1))\in H^2_{\mathrm{sing}}(X,\Z),
\]
and again denote by $\xi$ its reduction modulo $m$ in $H^2_{\mathrm{sing}}(X,A)$.
Restricting to a fibre $F\simeq \PP^1$, one has
\[
\xi|_F=c_1(\mathcal O_{\PP^1}(1))\neq 0\in H^2_{\mathrm{sing}}(\PP^1,A)\cong A.
\]
On the other hand, every class in $p^*H^2(C,A)$ restricts trivially to $F$, because
$p|_F:F\to C$ is constant. Hence
\[
\xi\notin p^*H^2(C,A)=H^2_{\mathrm{sing}}(X,A)_{\pi_1}.
\]
Therefore $c_X^*$ is not surjective, hence not an isomorphism. This proves \textup{(2)}.

We now prove \textup{(3)}. By the general necessary containment for semi-topological
degree-$2$ realizability,
\[
\mathrm{Im}(\Phi_X^2)\subseteq H^2_{\mathrm{sing}}(X,A)_{\pi_1}.
\tag{$**$}
\]
On the other hand, for the curve $C$ the degree-$2$ comparison map
\[
\Phi_C^2:H^2_{ST}(C,A)\longrightarrow H^2_{\mathrm{sing}}(C,A)
\]
is surjective. Let $\alpha\in H^2_{\mathrm{sing}}(C,A)$. Choose
$\widetilde{\alpha}\in H^2_{ST}(C,A)$ such that $\Phi_C^2(\widetilde{\alpha})=\alpha$.
By Theorem \ref{thm:comparison-naturality} and the functoriality of the comparison map with
respect to the morphism $p:X\to C$, the diagram
\[
\xymatrix{
H^2_{ST}(C,A) \ar[r]^{\Phi_C^2} \ar[d]_{p^*_{ST}} &
H^2_{\mathrm{sing}}(C,A) \ar[d]^{p^*} \\
H^2_{ST}(X,A) \ar[r]^{\Phi_X^2} &
H^2_{\mathrm{sing}}(X,A)
}
\]
commutes. Hence
\[
\Phi_X^2\bigl(p^*_{ST}(\widetilde{\alpha})\bigr)=p^*(\alpha).
\]
Therefore
\[
p^*H^2_{\mathrm{sing}}(C,A)\subseteq \mathrm{Im}(\Phi_X^2).
\tag{$***$}
\]
Combining $(*)$, $(**)$, and $(***)$, we obtain
\[
\mathrm{Im}(\Phi_X^2)
=
H^2_{\mathrm{sing}}(X,A)_{\pi_1}
=
p^*H^2_{\mathrm{sing}}(C,A),
\]
which proves \textup{(3)}.

Finally, let $c_0\in C$ be a point, and let $F:=p^{-1}(c_0)\subset X$
be the fibre divisor. Its first Chern class satisfies
\[
c_1(\mathcal O_X(F))=p^*c_1(\mathcal O_C(c_0)).
\]
Reducing modulo $m$, we see that the generator of $p^*H^2(C,A)\cong A$ lies in
$\rho_m(NS(X))$. Hence
\[
p^*H^2(C,A)\subseteq \rho_m(NS(X)).
\]
Using \textup{(3)}, we conclude that
\[
\rho_m(NS(X))\cap H^2_{\mathrm{sing}}(X,A)_{\pi_1}
=
p^*H^2(C,A)
=
\rho_m(NS(X))\cap \mathrm{Im}(\Phi_X^2).
\]
This shows that the $\pi_1$-detectable Weierstrass realizability conjecture holds.
\end{proof}



\begin{thebibliography}{99}

\bibitem{BerhuyGaloisCohomology}
G.~Berhuy,
\newblock {\em An Introduction to Galois Cohomology and its Applications},
\newblock London Mathematical Society Lecture Note Series, Vol.~377,
\newblock Cambridge University Press, Cambridge, 2010.

\bibitem{BorceuxJanelidzeGalois}
F.~Borceux and G.~Janelidze,
\newblock {\em Galois Theories},
\newblock Cambridge Studies in Advanced Mathematics, Vol.~72,
\newblock Cambridge University Press, Cambridge, 2001.



\bibitem{HansenBraidsCoverings}
V.~L.~Hansen,
\newblock \emph{Braids and Coverings: Selected Topics},
\newblock London Mathematical Society Student Texts, Vol.~18,
Cambridge University Press, 1989.


\bibitem{LT15}
H.-Y.~Liao and J.-H.~Teh,
\newblock ``Semi-topological Galois theory and the inverse Galois problem,''
\newblock Algebra Colloquium \textbf{22} (2015), no.~4, 687--706.

\bibitem{NeukirchSchmidtWingberg}
J.~Neukirch, A.~Schmidt and K.~Wingberg,
\newblock {\em Cohomology of Number Fields},
\newblock Grundlehren der mathematischen Wissenschaften, Vol.~323,
\newblock 2nd ed., Springer-Verlag, Berlin, 2008.



\bibitem{RahmanSchmeisser2002}
Q.~I.~Rahman and G.~Schmeisser,
\newblock {\em Analytic Theory of Polynomials},
\newblock London Mathematical Society Monographs. New Series, Vol.~26,
\newblock Oxford University Press, Oxford, 2002.

\bibitem{SerreLRFG}
J.-P.~Serre,
\newblock \emph{Linear Representations of Finite Groups},
\newblock Graduate Texts in Mathematics, Vol.~42,
Springer-Verlag, New York–Heidelberg, 1977.


\bibitem{SerreGaloisCohom}
J.-P.~Serre,
\newblock {\em Galois Cohomology},
\newblock Springer Monographs in Mathematics,
\newblock Springer-Verlag, Berlin, 1997.


\bibitem{RibesZalesskii}
L.~Ribes and P.~Zalesskii,
\newblock {\em Profinite Groups},
\newblock Ergebnisse der Mathematik und ihrer Grenzgebiete. 3. Folge. A Series of Modern Surveys in Mathematics, Vol.~40,
\newblock Springer-Verlag, Berlin, 2000.


\bibitem{SGA1}
A. Grothendieck,
\emph{Rev\^etements \'{e}tales et groupe fondamental (SGA 1)},
Lecture Notes in Math. 224, Springer, 1971.

\end{thebibliography}
\end{document}